\documentclass{article}

\usepackage{authblk}

\usepackage[utf8]{inputenc}
\usepackage[T1]{fontenc}
\usepackage[english]{babel}
\usepackage{textcomp}
\usepackage{amsmath,amssymb,amsthm}
\usepackage{lmodern}
\usepackage[a4paper]{geometry}
\usepackage{xcolor, pict2e}
\usepackage{microtype}
\usepackage{caption}
\usepackage{listings}
\usepackage{multicol}
\usepackage{moreverb}
\usepackage{hyperref}
\hypersetup{pdfstartview=XYZ}
\usepackage{wrapfig}
\usepackage[sans]{dsfont}
\usepackage{booktabs}

\usepackage{mathrsfs}
\usepackage{tikz, pgfplots}
\usepgfplotslibrary{groupplots}
\usetikzlibrary{positioning}
\usepackage{amsmath,amssymb}
\usetikzlibrary{decorations.pathreplacing}

\usepackage{ytableau}
\usepackage{array}

\usepackage{mathdots}
\usepackage{hyperref}
\hypersetup{
    colorlinks=true,
    linkcolor=blue,
    citecolor=magenta,
    urlcolor=blue,
    pdfborder={0 0 0}
}

\definecolor{fond}{rgb}{0.05,0.05,0.25}

\DeclareMathOperator{\ch}{\mathrm{ch}}
\DeclareMathOperator{\Id}{\mathrm{Id}}

\DeclareMathOperator{\cyc}{\mathrm{cyc}}

\DeclareMathOperator{\Pl}{\mathrm{Pl}}

\newcommand{\URS}{\ensuremath{\textup{URS}}}

\newcommand{\kS}{{\ensuremath{\mathfrak{S}}} }

\usepackage[font=sf, labelfont={sf,bf}, margin=1cm]{caption}
\usepackage{enumitem}
\setlist[itemize,1]{nosep}
\setlist[enumerate,1]{itemsep=0pt,label=(\alph*)}

\usepackage{float}
\usepackage[caption = false]{subfig}
\usepackage{graphicx}

\newtheorem*{theorem*}{Theorem}
\newtheorem{theorem}{Theorem}[section]
\newtheorem{proposition}[theorem]{Proposition}
\newtheorem{lemma}[theorem]{Lemma}
\newtheorem{corollary}[theorem]{Corollary}
\newtheorem{example}[theorem]{Example}

\theoremstyle{remark}
\newtheorem{remark}{Remark}[section]

\usepackage{todonotes}

\newcommand{\cE}{{\ensuremath{\mathcal E}} }

\newcommand{\bbE}{{\ensuremath{\mathbb E}} }

\newcommand{\bbZ}{{\ensuremath{\mathbb Z}} }

\newcommand{\bbR}{{\ensuremath{\mathbb R}} }
\newcommand{\bbC}{{\ensuremath{\mathbb C}} }

\newcommand{\ag}{\left\{ } 

\newcommand{\ad}{\right\} }

\newcommand{\cg}{\left[}
\newcommand{\cd}{\right]}
\newcommand{\pg}{\left(} 
\newcommand{\pd}{\right)}
\newcommand{\bg}{\left|}
\newcommand{\bd}{\right|}
\newcommand{\lf}{\left\lfloor}
\newcommand{\rf}{\right\rfloor}

\newcommand{\du}{{\ensuremath{\;:\;}}} 

\DeclareMathOperator{\ST}{ST}
\DeclareMathOperator{\RT}{RT}
\DeclareMathOperator{\height}{ht}
\DeclareMathOperator{\supp}{supp}

\makeatletter
\newcommand*\bigcdot{\mathpalette\bigcdot@{.5}}
\newcommand*\bigcdot@[2]{\mathbin{\vcenter{\hbox{\scalebox{#2}{$\m@th#1\bullet$}}}}}
\makeatother

\numberwithin{equation}{section}

\pgfplotsset{compat=1.18}

\begin{document}

\renewcommand{\theparagraph}{\thesubsection.\arabic{paragraph}} 
\title{Bounds on skew dimensions and characters of symmetric groups via thick hook decompositions}
\author{Lucas Teyssier}
\affil{Université de Lorraine, \texttt{lucas.teyssier@univ-lorraine.fr}}
\date{}

\maketitle

%

\begin{abstract}
  We introduce a new decomposition of Young diagrams into thick hooks. These thick hook decompositions enable having a better control over hook products for excited diagrams in the Naruse hook length formula, which leads to improved bounds on the number of standard tableaux of skew shapes. Combining this with elementary counting arguments in the Murnaghan--Nakayama rule, we establish a uniform bound on characters of symmetric groups $\mathfrak{S}_n$. In the case of balanced representations, this improves on the character bounds of Féray and Śniady for permutations with support size at least $n^{2/3}$, and is sharp for permutations with support size of order $n$. Finally, we recover some bounds of Pak and Panova on Kronecker coefficients with very short proofs that rely on characters.
\end{abstract}

\section{Introduction}
Let $\lambda$ be an irreducible representation of $\kS_n$ and $\sigma \in \kS_n$. Denote the dimension of $\lambda$ by $d_\lambda$, the character of $\lambda$ at $\sigma$ by $\ch^\lambda(\sigma)$, the associated renormalized character by $\chi^\lambda(\sigma) := \ch^\lambda(\sigma)/d_\lambda$, the minimal number of transpositions one needs to multiply to obtain $\sigma$ by $|\sigma|$, and the support size of $\sigma$ by $\supp(\sigma)$. 
The irreducible representations of $\kS_n$ are in bijection with integer partitions $\lambda \vdash n$, which can be represented by Young diagrams, and we use the same notation $\lambda$ for these three notions. The first row (resp. column) of a diagram $\lambda$ is denoted by $\lambda_1$ (resp. $\lambda'_1$), and the maximal hook length of $\lambda$ is $s(\lambda) := \lambda_1 + \lambda_1' - 1$.
Finally, for $C>0$ we say that $\lambda \vdash n$ is $C$-balanced if $s(\lambda) \leq C \sqrt{n}$.

\medskip

Building on the work of Kerov \cite{Kerov1998InterlacingMeasures}, Biane \cite{Biane1998RepresentationsFreeProbability,Biane2001ApproximateFactorizationConcentration} proved, alongside factorization formulas for characters, that for balanced diagrams $\lambda$ and permutations $\sigma$ with $\supp(\sigma) = O(1)$,
\begin{equation}
\label{eq: conjecture de Moore et Russell}
    |\chi^\lambda(\sigma)| \leq \pg \frac{O(1)}{\sqrt{n}}\pd^{|\sigma|}.
\end{equation}
This bound on characters for balanced diagrams was later extended to $\supp(\sigma) = O(\sqrt{n})$ by Rattan and Śniady \cite{RattanSniady2008}. Moore and Russell conjectured (\cite[Conjecture 3]{RattanSniady2008}) that \eqref{eq: conjecture de Moore et Russell} holds for any support size. Note that we always have $\supp(\sigma)/2 \leq |\sigma| \leq \supp(\sigma)-1$, so $\supp(\sigma)\asymp |\sigma|$.

The known bounds for $\supp(\sigma) = \omega(\sqrt{n})$ are however still relatively weak, especially for larger support sizes. Rattan and Śniady \cite{RattanSniady2008} proved that for permutations with $\supp(\sigma) = \omega(\sqrt{n})$, 
\begin{equation}\tag{RŚ}\label{eq: bound of Rattan Śniady for large support}
    |\chi^\lambda(\sigma)| \leq \pg O(1)\frac{|\sigma|^2}{n^{3/2}}\pd^{|\sigma|},
\end{equation}
which is trivial if $\supp(\sigma) = \omega(n^{3/4})$.

The strongest known character bound in this setting is due to Féray and Śniady \cite{FeraySniady2011}, who proved the following bound, which applies to all shapes of diagrams:
\begin{equation}\tag{FŚ}\label{eq: bound of Feray Śniady}
    |\chi^\lambda(\sigma)| \leq \pg O(1) \max\pg \frac{|\sigma|}{n}, \frac{s(\lambda)}{n}\pd\pd^{|\sigma|}.
\end{equation}
For balanced diagrams, \eqref{eq: bound of Feray Śniady} coincides with the bound of Rattan and Śniady if $|\sigma| = O(\sqrt{n})$; these bounds are excellent and led to optimal applications to the number of standard diagrams of balanced skew Young tableaux in \cite{DousseFeray2019skewdiagramscharacters}. If $|\sigma| = \omega(\sqrt{n})$, and still for balanced diagrams $\lambda$, \eqref{eq: bound of Feray Śniady}  can be rewritten as
\begin{equation}\tag{FŚbal}\label{eq: bound of Feray Śniady for large supports and balanced diagrams}
    |\chi^\lambda(\sigma)| \leq \pg O(1) \frac{|\sigma|}{n}\pd^{|\sigma|}.
\end{equation}
This bound is stronger than \eqref{eq: bound of Rattan Śniady for large support} since $|\sigma|/n = o(|\sigma|^{2}/n^{3/2})$.

Our main result is the following character bound.

\begin{theorem}\label{thm: borne avec racine de sigma cas général}
    As $n\to \infty$, uniformly over all representations $\lambda \vdash n$ and all permutations $\sigma \in \kS_n\backslash\ag \Id\ad$, we have
\begin{equation}
\label{eq: borne avec racine de sigma cas général}
     |\chi^\lambda(\sigma)| \leq \pg \frac{O(1)}{\sqrt{|\sigma|}}\pd^{|\sigma|} \max \pg  1 , \frac{s(\lambda) \sqrt{|\sigma|}}{n}\pd^{\supp(\sigma)}.
\end{equation}
\end{theorem}

Theorem \ref{thm: borne avec racine de sigma cas général} has the following simpler form for balanced representations.

\begin{theorem}\label{thm: borne avec racine de sigma cas équilibré}
     Let $C>0$. As $n\to \infty$, uniformly over all $C$-balanced representations $\lambda \vdash n$ and all permutations $\sigma \in \kS_n\backslash\ag \Id\ad$, we have
\begin{equation}
\label{eq: borne avec racine de sigma cas équilibré}
     |\chi^\lambda(\sigma)| \leq \pg \frac{O(1)}{\sqrt{|\sigma|}}\pd^{|\sigma|}.
\end{equation}
\end{theorem}

An interesting aspect of our results is that they are more precise for permutations with larger support sizes, contrary to the bounds of \cite{RattanSniady2008} and \cite{FeraySniady2011}. More precisely,
Theorem \ref{thm: borne avec racine de sigma cas équilibré} is stronger than \eqref{eq: bound of Feray Śniady for large supports and balanced diagrams} if $|\sigma| = \omega(n^{2/3})$, but weaker if $ |\sigma| =o(n^{2/3})$.
As such, Theorem \ref{thm: borne avec racine de sigma cas général} complements the Féray--Śniady character bounds, as we illustrate in Figure \ref{fig: character bounds for balanced diagrams all in one}.
Moreover, Theorem \ref{thm: borne avec racine de sigma cas équilibré} agrees with the conjecture of Moore and Russell for permutations with $\supp(\sigma) \asymp n$.

\medskip

Theorems \ref{thm: borne avec racine de sigma cas général} and \ref{thm: borne avec racine de sigma cas équilibré} are proved in Section \ref{s: proof main character bound}. To this end, we establish two intermediate results, which are of independent interest.
The first one is a bound on the number $d_{\lambda\backslash \mu}$ of standard tableaux of a skew diagram $\lambda\backslash\mu$, in function of the maximal hook length $s(\lambda)$ of $\lambda$ and the sizes of $\lambda$ and $\mu$. 

\begin{theorem}\label{thm: bound on skew dimensions intro}
     There exists a universal constant $C>0$ such that the following holds. Let $n\geq 1$ and $1\leq k\leq n$ be integers. Let $\lambda \vdash n$ and $\mu \vdash k$ such that $\mu \subset \lambda$. Then 
    \begin{equation}
         \frac{d_{\lambda\backslash \mu}}{d_\lambda} \leq  \pg C \max \pg  \frac{1}{\sqrt{k}} , \frac{s(\lambda)}{n}\pd \pd^k. 
    \end{equation} 
\end{theorem}
Theorem \ref{thm: bound on skew dimensions intro} is proved in Section \ref{s: bounds on skew dimensions}. It relies on thick hook decompositions in the Naruse hook length formula. We also show in Section \ref{s: sharpness of the bound} that Theorem \ref{thm: bound on skew dimensions intro} is sharp for all orders of magnitude for $k$ and $s(\lambda)$.

\medskip

Our second intermediate result is a simple character bound, written in terms of the number of rows $\ell(\lambda)$ -- or the diagonal length $\delta(\lambda)$ -- of a diagram $\lambda$ and the number of cycles $\cyc(\sigma)$ of a permutation $\sigma$.
\begin{theorem}\label{thm: character bound useful for fixed point free intro}
    Let $n\geq 1$, $\lambda \vdash n$, and $\sigma \in \kS_n$. Then 
\begin{enumerate}
    \item $|\ch^\lambda(\sigma)| \leq \ell(\lambda)^{\cyc(\sigma)}$;
    \item $|\ch^\lambda(\sigma)| \leq (2 \delta(\lambda))^{\cyc(\sigma)}$.
\end{enumerate}    
\end{theorem}
Theorem \ref{thm: character bound useful for fixed point free intro} is proved in Section \ref{s: Character bound in function of the diagonal length and the number of cycles}. Its proof relies on the Murnaghan--Nakayama rule and elementary counting arguments. We also use it in Section \ref{s: Kronecker} to obtain bounds on Kronecker coefficients.

\begin{figure}[ht]
\centering

\begin{minipage}[t]{0.49\textwidth}
\centering
\vspace{0pt}

\begin{tikzpicture}[x=4.5cm,y=4.5cm]

\draw[lightgray,densely dashed,line width=0.4pt]
    (0.5,0) -- (0.5,0.55);

\draw[lightgray,densely dashed,line width=0.4pt]
    (0.666666,0) -- (0.666666,0.55);

\draw[lightgray,densely dashed,line width=0.4pt]
    (0,0.333333) -- (1,0.333333);

\draw[lightgray,densely dashed,line width=0.4pt]
    (0,0.5) -- (1,0.5);

\draw[->] (0,0) -- (0,0.62);
\draw (0,0) -- (1,0);

\node[anchor=south west,font=\small] at (0,0.62)
    {$\displaystyle
      y=\frac{\ln(B^{-1/|\sigma|})}{\ln n}$};

\node[anchor=north east,font=\small] at (1.4,0.10)
    {$\displaystyle
      x=\frac{\ln|\sigma|}{\ln n}$};

\foreach \x/\lab in {
    0.5/{1/2},
    0.666666/{2/3},
    1/{1}
}{
    \draw (\x,-0.012) -- (\x,0.012);
    \node[anchor=north,font=\scriptsize]
        at (\x,-0.025) {$\lab$};
}

\foreach \y/\lab in {
    0.333333/{1/3},
    0.5/{1/2}
}{
    \draw (-0.010,\y) -- (0.010,\y);
    \node[anchor=east,font=\scriptsize]
        at (-0.020,\y) {$\lab$};
}

\node[anchor=north east,font=\small]
    at (-0.008,-0.012) {$0$};

\draw[
    very thick,
    pink
]
    (0,0.5)
    -- (0.5,0.5)
    -- (1,0);

\draw[
    very thick,
    blue!85!black
]
    (0,0)
    -- (1,0.5);

\draw[
    thick,
    densely dashed,
    black!60
]
    (0,0.5)
    -- (1,0.5);

\end{tikzpicture}

\end{minipage}
\hfill
\begin{minipage}[t]{0.49\textwidth}
\centering
\vspace{1.7cm}

\begin{tabular}{@{}cl@{}}

\tikz[baseline=-0.5ex]
    \draw[thick,densely dashed,black!60]
    (0,0)--(0.65,0);
&
Conjectured bound (Moore--Russell)
\\[5pt]

\tikz[baseline=-0.5ex]
    \draw[very thick,pink]
    (0,0)--(0.65,0);
&
Previous best bound (Féray--Śniady)
\\[5pt]

\tikz[baseline=-0.5ex]
    \draw[very thick,blue!85!black]
    (0,0)--(0.65,0);
&
Our bound from Theorem \ref{thm: borne avec racine de sigma cas équilibré}

\end{tabular}

\end{minipage}

\caption{Asymptotic uniform character bounds for balanced diagrams,
logarithmically rescaled. Here $B$ denotes the upper bound on
$\bg \chi^\lambda(\sigma)\bd$ from the right-hand-side of
\eqref{eq: conjecture de Moore et Russell} for the conjecture,
\eqref{eq: bound of Feray Śniady for large supports and balanced diagrams}
for Féray--Śniady, and
\eqref{eq: borne avec racine de sigma cas équilibré} for our bound.
The implicit universal constants (in the definition of balanced diagrams
and the character bounds) have a negligible asymptotic contribution,
with the rescaling of the figure.}
\label{fig: character bounds for balanced diagrams all in one}

\end{figure}

\section*{Notation}
 We use the following asymptotic notation. We write $f(n) = O(g(n))$ or $f(n) \lesssim g(n)$ or $g(n) = \Omega(f(n))$ if there exists a constant $C>0$ such that $|f(n)| \leq C |g(n)|$ for all $n$ large enough, $f(n) = o(g(n))$ or $g(n) = \omega(f(n))$
 if $f(n)/g(n) \xrightarrow[n\to \infty]{}0$, and $f(n) \asymp g(n)$ if $f(n) \lesssim g(n)$ and $f(n) \gtrsim g(n)$.
We use the French convention for Young diagrams. 

\section{Bounds on skew dimensions}\label{s: bounds on skew dimensions}
Let $\lambda$ and $\mu \subset \lambda$ be Young diagrams. We denote the number of standard tableaux of the skew diagram $\lambda\backslash \mu$ by $d_{\lambda\backslash \mu}$, and refer to this quantity as \textit{skew dimension}.
Different explicit formulas have been found to compute skew dimensions, see for instance \cite{OkunkovOlshanski1997ShiftedSchurFunctions, Stanley2003EnumerationSkewYoungTableaux}.

A decade ago, Naruse \cite{Naruse2014} found a way to compute skew dimensions via excited diagrams, a notion introduced in \cite{IkedaNaruse2009exciteddiagrams}. This powerful formula, now called Naruse hook length formula, has seen other proofs and generalizations, see \cite{MoralesPakPanova2018I, Konvalinka2020naruse}, and has enabled computing improved bounds on skew dimensions in \cite{MoralesPakPanova2018asymptotics,MoralesPakTassy2022skewshapelozenge, Olesker-TaylorTeyssierThevenin2025sharpboundsandcutoff} as well as explicit formulas for various skew shapes \cite{KimOh2017SelbergIntegralAndYoungBooks, MoralesPakPanova2019IIIProductFormulas, KimYoo2020ProductFormulasSkewTableaux}.
In this section we are concerned in proving general bounds on the quantity $\frac{d_{\lambda\backslash\mu}}{d_\lambda}$ for all shapes, in function of the sizes of the diagrams $\lambda$ and $\mu$ and the maximal hook length of $\lambda$.

\subsection{The Naruse hook length formula}

Let us now define excited diagrams and state the Naruse hook length formula. We represent boxes by the coordinates of their top-right corner in $\bbZ^2$.
If $u = (i,j) \in \bbZ^2$ is a box, we define its upper right square by $\URS(u) = \ag (i,j), (i+1, j), (i,j+1), (i+1,j+1) \ad$. Let $\lambda$ be a Young diagram, and let $E$ be a subset (of the boxes) of $\lambda$. A box $u\in E$ is said to be $(\lambda,E)$-\textit{excitable} if $\URS(u)\subset \lambda$ and $E\cap \URS(u) = \ag u \ad$. A \textit{simple excitation} (of $E$ in $\lambda$) consists of transforming the set of boxes $E$ into the diagram $E':= (E\backslash\ag u \ad)\cup \ag (i+1,j+1) \ad$, where  $u=(i,j)$ is an excitable box of $E$ in $\lambda$. In words, $E'$ has the same boxes as $E$ except one that was moved diagonally towards north-east. Finally an \textit{excited diagram} of a set of boxes $\mu$ in $\lambda$ is a diagram obtained by applying a (possibly empty) sequence of simple excitations. We denote the set of all excited diagrams of $\mu$ in $\lambda$ by $\cE(\lambda, \mu)$, and illustrate this in Figure \ref{fig: set of excited diagrams}.
We refer to \cite{Olesker-TaylorTeyssierThevenin2025sharpboundsandcutoff} or \cite{MoralesPakPanova2018I} for a more detailed description and further examples. 

\begin{figure}[!ht]
    \begin{center}
\ytableausetup{boxsize=1.2em}
\begin{ytableau}
    \none   &*(darkgray)& \\
    \none   &*(darkgray)&&&& \\
    \none   &*(darkgray)&*(lightgray)&&&\\
    \none   &*(darkgray)&*(darkgray)&*(lightgray)&&
\end{ytableau}
\begin{ytableau}
    \none   &*(darkgray)& \\
    \none   &*(darkgray)&&*(darkgray)&& \\
    \none   &*(darkgray)&&&&\\
    \none   &*(darkgray)&*(darkgray)&*(lightgray)&&
\end{ytableau}
\begin{ytableau}
    \none   &*(darkgray)& \\
    \none   &*(darkgray)&&&& \\
    \none   &*(darkgray)&*(lightgray)&&*(lightgray)&\\
    \none   &*(darkgray)&*(darkgray)&& &
\end{ytableau}
\begin{ytableau}
    \none   &*(darkgray)& \\
    \none   &*(darkgray)&&*(darkgray)&& \\
    \none   &*(darkgray)&&&*(lightgray)&\\
    \none   &*(darkgray)&*(lightgray)&& &
\end{ytableau}

\medskip

\begin{ytableau}
    \none   &*(darkgray)& \\
    \none   &*(darkgray)&&*(darkgray)&& \\
    \none   &*(darkgray)&&*(darkgray)&*(lightgray)&\\
    \none   &*(darkgray)&&& &
\end{ytableau}
\begin{ytableau}
    \none   &*(darkgray)& \\
    \none   &*(darkgray)&&&&*(darkgray) \\
    \none   &*(darkgray)&*(lightgray)&&&\\
    \none   &*(darkgray)&*(darkgray)&& &
\end{ytableau}
\begin{ytableau}
    \none   &*(darkgray)& \\
    \none   &*(darkgray)&&*(darkgray)&&*(darkgray) \\
    \none   &*(darkgray)&&&&\\
    \none   &*(darkgray)&*(lightgray)&& &
\end{ytableau}
\begin{ytableau}
    \none   &*(darkgray)& \\
    \none   &*(darkgray)&&*(darkgray)&&*(darkgray) \\
    \none   &*(darkgray)&&*(darkgray)&&\\
    \none   &*(darkgray)&&& &
\end{ytableau}
\end{center}
\caption{The set $\cE(\lambda, \mu)$, with $\lambda = [5,5,5,2]$ and $\mu = [3,2,1,1]$, consists of the following diagrams. We represent the excitable boxes light gray and the other boxes of excited diagrams in dark gray.}
\label{fig: set of excited diagrams}
\end{figure}

\medskip

We denote the hook length of a box $u$ in $\lambda$ by $H(\lambda, u)$. If $E\subset \lambda$ is a subset of the boxes of $\lambda$ we set $H(\lambda, E) = \prod_{u\in E}H(\lambda, u)$. In particular, the hook length formula rewrites as $d_\lambda = |\lambda|!/H(\lambda, \lambda)$.
We use the falling factorial notation: if $0\leq k \leq n$ are two integers, we write 
\begin{equation}
    n^{\downarrow k} = n(n-1)\cdots(n-k+1) = \frac{n!}{(n-k)!}.
\end{equation}
Finally, let $\mu\subset \lambda$ be two Young diagrams. We define  the \textit{excited sum} of $\mu$ in $\lambda$ by
\begin{equation}
    S(\lambda, \mu) = \sum_{E\in \cE(\lambda, \mu)} H(\lambda, E).
\end{equation}
The Naruse hook length formula can then be written as the following identity: 
\begin{equation}\label{eq: NHLF forme avec somme excitée}
    \frac{d_{\lambda\backslash \mu}}{d_\lambda} = \frac{S(\lambda, \mu)}{|\lambda|^{\downarrow |\mu|}}.
\end{equation}

\subsection{Motivation towards decompositions of diagrams}

A first way to bound excited sums $S(\lambda, \mu)$ is to bound each term of the sum by the maximal summand. Then we obtain the bound

\begin{equation}\label{eq: max sum bound for excited sums}
    S(\lambda, \mu) \leq |\cE(\lambda, \mu)| H(\lambda, \mu).
\end{equation}

Let us illustrate the (non-)sharpness of this bound when $\mu = [1]$ consists of only one box. Then have $d_{\lambda\backslash \mu} = d_{\lambda}$ and (equivalently) $S(\lambda, \mu) = n$. For example, the identity $S(\lambda, \mu) = n$ follows from the fact that the union of the hooks of the boxes with coordinates $(1,1), (2,2), \ldots$ in $\lambda$ is exactly $\lambda$.

\begin{example}\label{ex: illustration non sharness of max sum}
    In the following 3 cases, we omit integer parts, and the equivalents are as $n\to \infty$.
 \begin{enumerate}
     \item Balanced shape, represented here by the square diagram. If $\lambda = [\sqrt{n}^{\sqrt{n}}]$, then $|\cE(\lambda, \mu)| = \sqrt{n}$ and $H(\lambda, \mu) = 2 \sqrt{n} - 1$, so \eqref{eq: max sum bound for excited sums} gives $S(\lambda, \mu)\leq 2n$, which is rather good.
     \item Thin shape, represented here by a diagram with only two rows, of the same length. If $\lambda = [n/2, n/2]$, then $|\cE(\lambda, \mu)| = 2$ and $H(\lambda, \mu) = n/2 + 1$, so \eqref{eq: max sum bound for excited sums} gives $S(\lambda, \mu)\leq n(1+o(1))$, which is very good.
     \item Hybrid shape: union of a long row and a square. If $\lambda = [3n/4, m^m]$, where $m:= \sqrt{n/4}$, then $|\cE(\lambda, \mu)| = (1+o(1))\sqrt{n}/2$ and $H(\lambda, \mu) = (1+o(1))3n/4$, so \eqref{eq: max sum bound for excited sums} gives $S(\lambda, \mu)\leq (3/8+o(1))n^{3/2}$, which is off by a factor of order $n^{1/2}$.
     \end{enumerate}
\end{example}

This illustrates that \eqref{eq: max sum bound for excited sums} is good for diagrams whose rows have comparable length, but can be wasteful for diagrams with hybrid shapes.

In Example \ref{ex: illustration non sharness of max sum} (c), it is not hard to modify the argument to obtain a precise estimate: we can split the excited diagrams depending on whether the excited box is in position $(1,1)$ or in position $(i,i)$ for some $i\geq 2$. Writing again $m=\sqrt{n/4}$, we obtain this way the very good estimate
\begin{equation}
    S(\lambda, \mu) = S(\lambda^1, [1]) + S([m^m], [1]) = (3/4 + o(1))n + m^2 = (1+ o(1))n.
\end{equation}
Now we want to generalize this idea to non-trivial cases, where $\lambda$ and $\mu$ can have any size and shape. We will decompose the diagram $\lambda$ into parts that have comparable contributions to $S(\lambda, \mu)$.

\subsection{Thick hook decomposition of Young diagrams}

Recall that we denote the hook of the diagonal box $(i,i)$ in $\lambda$ by $\lambda^i$. For $j\leq k$ we set $\lambda^{j\to k} = \cup_{j\leq i\leq k} \lambda^i$. We illustrate this in Figure \ref{fig: thick hook definition}.

\begin{figure}[!ht]
    \begin{center}
        \begin{ytableau}
            \none & \\
            \none && \\
            \none &&&&&*(gray) \\
            \none &&&&&*(gray)&& \\
            \none &&&&&*(gray)&*(gray)&*(gray)&*(gray) \\
            \none &&&&&&&& \\
            \none &&&&&&&& \\
            \none &&&&&&&&&& \\
            \none &&&&&&&&&&& \\
        \end{ytableau}
        \begin{ytableau}
            \none & \\
            \none &&*(gray) \\
            \none &&*(gray)&*(gray)&*(gray)& \\
            \none &&*(gray)&*(gray)&*(gray)&&& \\
            \none &&*(gray)&*(gray)&*(gray)&&&& \\
            \none &&*(gray)&*(gray)&*(gray)&*(gray)&*(gray)&*(gray)&*(gray) \\
            \none &&*(gray)&*(gray)&*(gray)&*(gray)&*(gray)&*(gray)&*(gray) \\
            \none &&*(gray)&*(gray)&*(gray)&*(gray)&*(gray)&*(gray)&*(gray)&*(gray)&*(gray) \\
            \none &&&&&&&&&&& \\
        \end{ytableau}
    \end{center}
\caption{The diagram $\lambda = [11,10, 8^3, 7, 5, 2,1]$. The subset $\lambda^5$ is represented on the left and $\lambda^{2\to 4}$ is represented on the right, both in gray.}
\label{fig: thick hook definition}
\end{figure}

\subsubsection{Definition of thick hook decompositions}

Let $\lambda$ be a Young diagram, and $1\leq a\leq b$. An $(a,b)$ thick hook decomposition of $\lambda$ is a sequence of thick hooks $(T_j)_{1\leq j \leq p} = (\lambda^{(i_{j-1}+1) \to i_{j}})_{1\leq j \leq p}$, where $0= i_0<i_1<\ldots<i_p = \delta(\lambda)$, such that for every $1\leq j \leq p $ we have
    \begin{equation}
        a \leq |T_j| \leq b.
    \end{equation}
By definition $\lambda$ is the disjoint union of the thick hooks $T_j$. We illustrate this in Figure \ref{fig: thick hook decomposition}.

\begin{figure}[!ht]
     \begin{center}
        \begin{ytableau}
            \none &*(darkgray)\\
            \none &*(darkgray)&*(darkgray!70) &*(darkgray!40) \\
            \none &*(darkgray)&*(darkgray!70) &*(darkgray!40) &*(darkgray!40)&*(darkgray!10)&*(darkgray!10)\\
            \none &*(darkgray)&*(darkgray!70) &*(darkgray!40)&*(darkgray!40)&*(darkgray!10)&*(darkgray!10)&*(darkgray!10)\\
            \none &*(darkgray)&*(darkgray!70) &*(darkgray!40)&*(darkgray!40)&*(darkgray!10)&*(darkgray!10)&*(darkgray!10)&*(darkgray!10)&*(darkgray!10)\\
            \none &*(darkgray)&*(darkgray!70) &*(darkgray!40)&*(darkgray!40)&*(darkgray!10)&*(darkgray!10)&*(darkgray!10)&*(darkgray!10)&*(darkgray!10)&*(darkgray!10)\\
            \none &*(darkgray)&*(darkgray!70) &*(darkgray!40)&*(darkgray!40)&*(darkgray!10)&*(darkgray!10)&*(darkgray!10)&*(darkgray!10)&*(darkgray!10)&*(darkgray!10)&*(darkgray!10)\\
            \none &*(darkgray)&*(darkgray!70) &*(darkgray!40)&*(darkgray!40)&*(darkgray!40)&*(darkgray!40)&*(darkgray!40)&*(darkgray!40)&*(darkgray!40)&*(darkgray!40)&*(darkgray!40)&*(darkgray!40)\\
            \none &*(darkgray)&*(darkgray!70) &*(darkgray!40)&*(darkgray!40)&*(darkgray!40)&*(darkgray!40)&*(darkgray!40)&*(darkgray!40)&*(darkgray!40)&*(darkgray!40)&*(darkgray!40)&*(darkgray!40)&*(darkgray!40)&*(darkgray!40)\\
            \none &*(darkgray)&*(darkgray!70) &*(darkgray!70) &*(darkgray!70) &*(darkgray!70) &*(darkgray!70) &*(darkgray!70) &*(darkgray!70) &*(darkgray!70) &*(darkgray!70) &*(darkgray!70) &*(darkgray!70) &*(darkgray!70) &*(darkgray!70) &*(darkgray!70) &*(darkgray!70) &*(darkgray!70) &*(darkgray!70) &*(darkgray!70) \\
            \none &*(darkgray)&*(darkgray)&*(darkgray)&*(darkgray)&*(darkgray)&*(darkgray)&*(darkgray)&*(darkgray)&*(darkgray)&*(darkgray)&*(darkgray)&*(darkgray)&*(darkgray)&*(darkgray)&*(darkgray)&*(darkgray)&*(darkgray)&*(darkgray)&*(darkgray)&*(darkgray)&*(darkgray)&*(darkgray)&*(darkgray)&*(darkgray)\\
        \end{ytableau}
    \end{center}
\caption{A $(18, 36)$ thick hook decomposition of the diagram $\lambda = [24,19, 14, 12, 11, 10, 9, 7, 6, 3, 1]$. Here we have $(i_1, i_2, i_3,i_4) = (1,2, 4, 7)$ and $(|T_j|)_{1\leq j \leq 4} = (34, 26, 33, 23)$. In this case this is the unique such decomposition.}
\label{fig: thick hook decomposition}
\end{figure}

\subsection{Minimally excited diagrams in a thick hook decomposition}

Let $\lambda$ be a Young diagram. Let $1\leq a\leq b$ such that there is at least one thick hook decomposition of $\lambda$, and let $\boldsymbol{T} = (T_j)_{1\leq j \leq p}$ be such a decomposition. Let $1\leq \ell \leq \lambda_1$ be an integer, $\boldsymbol{\ell} = (\ell_j)_{1\leq j \leq p}$ be a sequence of non-negative integers such that $\sum_{j=1}^p \ell_j = \ell$, and let
\begin{equation}
    \cE(\lambda ; \boldsymbol{T}, \boldsymbol{\ell}) = \ag E \in \cE(\lambda, [\ell]) \mid E\cap T_j = \ell_j \text{ for all } 1\leq j \leq p \ad
\end{equation}
be the set of excited diagrams of $[\ell]$ in $\lambda$ with $\ell_j$ boxes in the thick hook $T_j$, for each $1\leq j\leq p$.

If $\cE(\lambda ; \boldsymbol{T}, \boldsymbol{\ell})$ is non-empty, all its diagrams are excited diagrams of a minimally excited diagram which is a union of shifted rows. 
Indeed, let $(i_j)_{0\leq j \leq p}$ be such that $T_j = \lambda^{(i_{j-1} + 1)\to i_j}$ for $1\leq i\leq p$, and let $m = \max \ag q\geq 1 \mid (\ell+q-1, q) \in \lambda \ad$. By the bijection with flagged tableaux from \cite[Section 3.3] {MoralesPakPanova2018I}, the diagrams in $\cE(\lambda ; \boldsymbol{T}, \boldsymbol{\ell})$ correspond exactly to sequences of integers $(u_1, \ldots, u_\ell)$ such that $1\leq u_1\leq \ldots \leq u_\ell \leq m$, and such that [$u_k \geq i_{j-1}+1$ whenever $\sum_{i=1}^{j-1} \ell_i < k \leq \sum_{i=1}^j \ell_i$ for some $1\leq j\leq p$], and the coordinate-wise partial order on vectors of $\bbR^\ell$ corresponds to the partial order on excited diagrams (where $E$ is smaller than $E'$ if $E'$ can be obtained from $E$ after applying excited moves). Among such sequences of integers, the sequence $(u_1, \ldots, u_\ell)$ defined by $u_k = i_{j-1}+1$ whenever $\sum_{i=1}^j \ell_j < k \leq \sum_{i=1}^j \ell_j$ for some $1\leq j\leq p$ is clearly the unique minimizer.

We denote this minimally excited diagram by $E_{\boldsymbol{T}, \boldsymbol{\ell}}$, and illustrate this in Figure \ref{fig: minimally excited diagram with respect to a thick hook decomposition}.

\begin{figure}[!ht]
    \begin{center}
        \begin{ytableau}
            \none &*(darkgray)\\
            \none &*(darkgray)&*(darkgray!70) &*(darkgray!40) \\
            \none &*(darkgray)&*(darkgray!70) &*(darkgray!40) &*(darkgray!40)&*(darkgray!10)&*(darkgray!10)\\
            \none &*(darkgray)&*(darkgray!70) &*(darkgray!40)&*(darkgray!40)&*(darkgray!10)&*(darkgray!10)&*(darkgray!10)\\
            \none &*(darkgray)&*(darkgray!70) &*(darkgray!40)&*(darkgray!40)&*(darkgray!10)&*(darkgray!10)&*(darkgray!10)&*(darkgray!10)&*(darkgray!10)\\
            \none &*(darkgray)&*(darkgray!70) &*(darkgray!40)&*(darkgray!40)&*(darkgray!10)&*(darkgray!10)&*(darkgray!10)&*(darkgray!10)&*(darkgray!10)&*(darkgray!10)\\
            \none &*(darkgray)&*(darkgray!70) &*(darkgray!40)&*(darkgray!40)&*(darkgray!10)&*(darkgray!10)&*(darkgray!10)&*(darkgray!10)&*(black)&*(darkgray!10)&*(darkgray!10)\\
            \none &*(darkgray)&*(darkgray!70) &*(darkgray!40)&*(darkgray!40)&*(darkgray!40)&*(darkgray!40)&*(darkgray!40)&*(darkgray!40)&*(darkgray!40)&*(darkgray!40)&*(darkgray!40)&*(darkgray!40)\\
            \none &*(darkgray)&*(darkgray!70) &*(darkgray!40)&*(black)&*(black)&*(black)&*(darkgray!40)&*(darkgray!40)&*(darkgray!40)&*(darkgray!40)&*(darkgray!40)&*(darkgray!40)&*(darkgray!40)&*(darkgray!40)\\
            \none &*(darkgray)&*(darkgray!70) &*(darkgray!70) &*(darkgray!70) &*(darkgray!70) &*(darkgray!70) &*(darkgray!70) &*(darkgray!70) &*(darkgray!70) &*(darkgray!70) &*(darkgray!70) &*(darkgray!70) &*(darkgray!70) &*(darkgray!70) &*(darkgray!70) &*(darkgray!70) &*(darkgray!70) &*(darkgray!70) &*(darkgray!70) \\
            \none &*(black)&*(darkgray)&*(darkgray)&*(darkgray)&*(darkgray)&*(darkgray)&*(darkgray)&*(darkgray)&*(darkgray)&*(darkgray)&*(darkgray)&*(darkgray)&*(darkgray)&*(darkgray)&*(darkgray)&*(darkgray)&*(darkgray)&*(darkgray)&*(darkgray)&*(darkgray)&*(darkgray)&*(darkgray)&*(darkgray)&*(darkgray)\\
        \end{ytableau}
    \end{center}
\caption{The minimally excited diagram of $\cE(\lambda ; \boldsymbol{T}, \boldsymbol{\ell})$, $E_{\boldsymbol{T}, \boldsymbol{\ell}}$, represented in black, where $\lambda = [24,19, 14, 12, 11, 10, 9, 7, 6, 3, 1]$, $\boldsymbol{T} = (\lambda^{1\to 1}, \lambda^{2\to 2}, \lambda^{3\to 4}, \lambda^{5\to 7})$, $\ell = 5$, and $\boldsymbol{\ell} = (1,0, 3, 1)$.}
\label{fig: minimally excited diagram with respect to a thick hook decomposition}
\end{figure}

\subsection{A general bound for line-shaped diagrams}
We first give a sufficient condition for thick hook decompositions to exist.
\begin{lemma}
    Let $\lambda$ be a Young diagram and $a\geq s := s(\lambda)$. Then there exists an $(a,4a)$ thick hook decomposition of $\lambda$.
\end{lemma}
\begin{proof}
    First build a $(1,4a)$ thick hook decomposition $\boldsymbol{T}$ of $\lambda$ by maximizing the size of $T_1$, then $T_2$, and so on up to a last thick hook $T_p$. Since $s\leq a$, by construction we have $3a \leq |T_j| \leq 4a$ for all $1\leq j \leq p-1$, and $|T_p| \leq 4a$. If $|T_p| \geq a$, then $\boldsymbol{T}$ is already a $(a,4a)$ thick hook decomposition of $\lambda$. Otherwise, using again that $s\leq a$, we can decrease the value of $i_{p-1}$ to transfer between $a$ and $2a$ boxes from $T_{p-1}$ to $T_p$, and after doing that we have $a\leq |T_{p-1}| \leq 3a$ and also $a\leq |T_p|\leq 3a$. This concludes the proof.
\end{proof}

We now extend a counting argument from \cite[Section 3.1]{Olesker-TaylorTeyssierThevenin2025sharpboundsandcutoff} to thick hook decompositions.

\begin{lemma}
    Let $\lambda$ be a Young diagram.  Let $a\geq s:=s(\lambda)$, where we recall that $s(\lambda)$ is the maximal hook length in $\lambda$. Let $\boldsymbol{T} = (T_j)_{1\leq j\leq p}$ be an $(a, 4a)$ thick hook decomposition of $\lambda$. Let $1\leq \ell \leq \lambda_1$. Then there are at most 
        $\binom{\ell + \lf n/a \rf}{\ell}$
    sequences of non-negative integers $\boldsymbol{\ell} = (\ell_j)_{1\leq j\leq p}$ such that $\sum_{j=1}^p \ell_j = \ell$ and $\cE(\lambda ; \boldsymbol{T}, \boldsymbol{\ell})$ is non-empty.
\end{lemma}
\begin{proof}
    Each $\boldsymbol{\ell}$ such that $\cE(\lambda ; \boldsymbol{T}, \boldsymbol{\ell})$ is non-empty corresponds to exactly one minimally excited diagram $E_{\boldsymbol{T}, \boldsymbol{\ell}}$, which is determined by the successive lengths $\ell_j$ and how many thick hooks they jump over. Moreover, we can represent each $E_{\boldsymbol{T}, \boldsymbol{\ell}}$ denoting successive boxes by $*$ and the jumps by $|$. (For example in Figure \ref{fig: minimally excited diagram with respect to a thick hook decomposition}, the configuration of black boxes can be represented by $*||***|*$.) Therefore since we have $\ell$ stars and $p-1$ bars, the number of such $\boldsymbol{\ell}$ is at most $\binom{\ell + (p-1)}{\ell}$, which is less than $\binom{\ell + p}{\ell}$.
     Finally, each $T_j$ contains at least $a$ boxes, and $p$ is an integer, so we have $p\leq \lf n/a \rf$. This concludes the proof.
\end{proof}

\begin{proposition}\label{prop: borne générale sur les sommes excitées avec a}
Let $\lambda$ be a Young diagram, and denote $n = |\lambda|$.  Let $a \in [s(\lambda), n]$, where we recall that $s(\lambda)$ is the maximal hook length in $\lambda$.  Let $1\leq \ell \leq \lambda_1$. Then 
\begin{equation}
    S(\lambda, [\ell]) \leq \binom{\ell + \lf n/a \rf }{\ell} (4a)^{\ell}.
\end{equation}
\end{proposition}
\begin{proof}
    Let $\boldsymbol{T} = (T_j)_{1\leq j\leq p}$ be an $(a, 4a)$ thick hook decomposition of $\lambda$.
Let $\boldsymbol{\ell} = (\ell_j)_{1\leq j \leq p}$ such that $\cE(\lambda ; \boldsymbol{T}, \boldsymbol{\ell})$ is non-empty.
Let $1\leq j\leq p$. Denote  the part of $E_{\boldsymbol{T}, \boldsymbol{\ell}}$ which is in $T_j$ by $\nu_j$.
For any set of boxes $\rho\subset T_j$, extending the notion of excited diagrams to any set of boxes, we set
\begin{equation}
    S^{*j}(\lambda, \rho) = \sum_E H(\lambda, E), 
\end{equation}
where the sum is over all excited diagrams of $\rho$ whose boxes are all in $T_j$.
Denote the bottom-left corner of $T_j$ by $u_0(j)$.
Then for each $1\leq j\leq p$ we have,
\begin{equation}
     S^{*j}(\lambda, \nu_j) \leq \prod_{u\in \nu_j} S^{*j}(\lambda, u) \leq \prod_{u\in \nu_j} S^{*j}(\lambda, u_0(j)) = \prod_{u\in \nu_j} |T_j| = |T_j|^{\ell_j} \leq (4a)^{\ell_j},
\end{equation}
and therefore
\begin{equation}
    \prod_{1\leq j \leq p} S^{*j}(\lambda, \nu_j) \leq  \prod_{1\leq j \leq p} (4a)^{\ell_j} = (4a)^{\ell}.
\end{equation}
We conclude, denoting the empty set of boxes by $\emptyset$, that
\begin{equation}
     S(\lambda, [\ell]) = \sum_{E\in \cE(\lambda, [\ell])} H(\lambda, E) = \sum_{\underset{\cE(\lambda ; \boldsymbol{T}, \boldsymbol{\ell}) \ne \emptyset}{\boldsymbol{\ell}= (\ell_j)_{1\leq j \leq p} \du }} \prod_{1\leq j \leq p} S^{*j}(\lambda, \nu_j) \leq \binom{\ell + \lf n/a \rf }{\ell} (4a)^{\ell}. \qedhere
\end{equation}
\end{proof}

\subsection{Specific bounds for line-shaped diagrams}

\begin{proposition}\label{prop: bound excited sums when mu is a line}
   Let $n\geq 1$ be an integer.  Let $\lambda \vdash n$. Recall that $s:= s(\lambda)$ denotes the maximal hook length in $\lambda$. Let $1\leq \ell \leq \lambda_1$.
\begin{enumerate}
    \item Assume that $\ell \leq n/s$. Then 
    \begin{equation}
        S(\lambda, [\ell]) \leq \pg 8\frac{n}{\ell}\pd^\ell.
    \end{equation}
    \item Assume that  $ \ell \geq n/s$. Then 
    \begin{equation}
        S(\lambda, [\ell]) \leq \pg 4e^{2}\pd^\ell s^\ell.
    \end{equation}
\end{enumerate}
\end{proposition}
\begin{proof}
    \begin{enumerate}
        \item Set $a = n/\ell$. By assumption we have $a\in [s, n]$. Moreover $\lf n/a \rf = \lf \ell\rf  = \ell$, so $\binom{\ell + \lf n/a \rf }{\ell} = \binom{2\ell}{\ell} \leq 2^\ell$ and we conclude using Proposition \ref{prop: borne générale sur les sommes excitées avec a} that 
        \begin{equation}
             S(\lambda, [\ell]) \leq \binom{\ell + \lf n/a \rf }{\ell}(4a)^\ell \leq (8a)^\ell = (8n/\ell)^\ell.
        \end{equation}
        \item Set $a = s$ and $q= \lf n/a\rf$. In particular we have $\lf n/a\rf\geq 2$ so $q\geq (n/a)/2$. Since we also have $q \leq \ell$ by assumption, it follows, writing $\gamma := q/\ell$, that
        \begin{equation}
            \binom{\ell + \lf n/a \rf }{\ell} = \binom{q + \ell}{q}  \leq \binom{2\ell}{q}= \frac{(2\ell)^{\downarrow q}}{q!} \leq \frac{(2\ell)^q}{(q/e)^q} = \pg\frac{2e\ell}{q}\pd^q = \pg\pg \frac{2e}{\gamma} \pd^{\gamma}\pd^\ell. 
        \end{equation}
    \end{enumerate}
But $(2e/\gamma)^{\gamma} \leq e^{2}$ since for any $A>0$ the maximum of the function $x\in \bbR^*_+ \mapsto (A/x)^x$ is $e^{A/e}$ (reached at $x = A/e$). It follows that
\begin{equation}
   \binom{\ell + \lf n/a \rf }{\ell} \leq e^{2\ell},
\end{equation}
and we conclude using Proposition \ref{prop: borne générale sur les sommes excitées avec a} that
\begin{equation}
    S(\lambda, [\ell]) \leq \binom{\ell + \lf n/a \rf }{\ell} (4a)^\ell \leq \pg 4e^2\pd^\ell s^\ell. \qedhere
\end{equation}
\end{proof}

\begin{remark}
To give some intuition on Proposition \ref{prop: bound excited sums when mu is a line}, in (a) $\ell$ is small so what dominates is the fact that boxes block one another, which explains why in the end we obtain about $1/\ell!$ times $n^\ell$, while in (b) there are very few excited diagrams (after the thick hook reduction), and what matters the most is the first term, $H(\lambda, [\ell])$, which is equal to $s^\ell$ up to an $e^{O(\ell)}$ error term.
\end{remark}

\begin{corollary}
   There exists a universal constant $C>0$ such that the following holds. Let $n\geq 1$ be an integer.  Let $\lambda \vdash n$, and denote again the maximal hook length in $\lambda$ by $s:= s(\lambda)$. Let $1\leq \ell \leq \lambda_1$. 
    \begin{equation}
        \frac{d_{\lambda\backslash[\ell]}}{d_\lambda}  \leq \pg C  \max\pg \frac{1}{\ell}, \frac{s}{n} \pd\pd ^\ell.
    \end{equation}
\end{corollary}
\begin{proof}
    This follows immediately from Proposition \ref{prop: bound excited sums when mu is a line} since $n^\ell/n^{\downarrow \ell} \leq \ell^\ell/\ell! \leq e^\ell$ for $1\leq \ell \leq n$.
\end{proof}

\begin{remark}
    We only proved bounds for row shapes $[\ell]$. By symmetry (up to taking conjugate diagrams, which doesn’t affect the maximal hook length), the statements above also apply to columns $[1^\ell]$.
\end{remark}

\subsection{Stairs-decomposition of a diagram}

Bounding the number of excited diagrams can be very delicate for general shapes. As we already observed, it is simpler when $\mu$ is line-shaped (i.e.\ $\mu =[\ell]$ or $\mu =[1^\ell]$ for some $\ell$). Then one can write the diagram $\mu$ as a union of not-too-many lines and deal with these lines separately.

As observed in \cite{Olesker-TaylorTeyssierThevenin2025sharpboundsandcutoff}, it is convenient to decompose a diagram $\mu$ as a union of at most $2\delta(\mu)$ lines, corresponding to the Frobenius coordinates of the diagram, but where we include also diagonal boxes to the rows. We call such a decomposition a stairs-decomposition, or $S$-decomposition. We illustrate this in Figure \ref{fig:stairs slicing}.

\begin{figure}[!ht]
    \begin{center}
    \begin{ytableau}
        \none &\none \\
        \none &*(darkgray!40) &\none &\none\\
        \none &*(darkgray!40)&*(darkgray!20)\\
        \none &*(darkgray!40)&*(darkgray!20)\\
        \none &*(darkgray!40)&*(darkgray!20)&*(darkgray!60)&*(darkgray!60)&*(darkgray!60)&*(darkgray!60)& \none &\none \\
        \none &*(darkgray!40)&*(darkgray!80)&*(darkgray!80)&*(darkgray!80)&*(darkgray!80)&*(darkgray!80)&*(darkgray!80)&*(darkgray!80) & \none &\none\\
        \none &*(darkgray)&*(darkgray)&*(darkgray)&*(darkgray)&*(darkgray)&*(darkgray)&*(darkgray)&*(darkgray)&*(darkgray)&*(darkgray)&*(darkgray)&*(darkgray)&*(darkgray)&*(darkgray) & \none &\none 
    \end{ytableau}
\end{center}

\caption{The stairs decomposition of $\lambda = [14,8,6,2,2,1]$.}
\label{fig:stairs slicing}
\end{figure}

\subsection{Bounds on skew dimensions for general shapes}

\begin{proposition}\label{prop: borne sommes excitées forme générale}
     There exists a universal constant $C>0$ such that the following holds. Let $1\leq k\leq n$ be integers. Let $\lambda \vdash n$ and $\mu \vdash k$ such that $\mu \subset \lambda$. Then 
    \begin{equation}
         S(\lambda, \mu) \leq  \pg C \max \pg \frac{n}{\sqrt{k}}, s \pd \pd^k. 
    \end{equation} 
\end{proposition}
\begin{proof}
    Denote by $(\ell_i)_{1\leq i\leq q}$ the parts in the stairs decomposition of $\mu$,  where we recall that the number of parts satisfies $q \leq 2\delta(\mu) \leq 2\sqrt{k}$. In this proof we use the notation $\ell_i$ for both the subsets of boxes in the stairs slicing and the lengths of these parts.

Let $I = \ag i\in [q] \mid \ell_i \leq n/s \ad$ and $J = [q]\backslash I$, where $[q] = \ag 1, 2, \ldots, q\ad$. Set $\ell_I = \sum_{i\in I} \ell_i$ and $\ell_J = \sum_{i\in J} \ell_i$. Note that $k = \ell_I + \ell_J$. Set also $C = \max\pg 8, 4e^2\pd = 4e^{2}$. By Proposition \ref{prop: bound excited sums when mu is a line}, we can reduce the problem to line-shaped diagrams, which leads to
\begin{equation}
     S(\lambda, \mu) \leq \prod_i S(\lambda, \ell_i) \leq \pg \prod_{i\in I} \pg C \frac{n}{\ell_i} \pd^{\ell_i}  \pd \pg \prod_{i\in J} \pg C s \pd^{\ell_i}  \pd = C^k n^{\ell_I} s^{\ell_J} \prod_{i\in I} \pg \frac{1}{\ell_i}\pd^{\ell_i}.
\end{equation}
Without loss of generality we may assume that $|I| \geq 1$. By convexity we have
\begin{equation}
    \prod_{i\in I} \pg \frac{1}{\ell_i}\pd^{\ell_i} \leq \prod_{i\in I} \pg \frac{1}{\ell_I/|I|}\pd^{\ell_I/|I|} = \pg \frac{|I|}{\ell_I}\pd^{\ell_I}.
\end{equation}
Therefore, since $|I|\leq q \leq 2\sqrt{k}$, 
\begin{equation}\label{eq: borne sommes excitées forme générale inter}
    S(\lambda, \mu) (2C)^{-k} \leq s^{\ell_J} \pg \frac{n|I|}{2\ell_I}\pd^{\ell_I} = s^k \pg \frac{n|I|}{2s\ell_I}\pd^{\ell_I} \leq s^k \pg \frac{n\sqrt{k}}{s\ell_I}\pd^{\ell_I} = s^k \pg \frac{n}{s\sqrt{k}} \frac{k}{\ell_I}\pd^{\ell_I}
\end{equation}
Now there are two cases. 
\begin{enumerate}
    \item If $\frac{n}{s\sqrt{k}} \leq e$, then by \eqref{eq: borne sommes excitées forme générale inter} and since $(A/x)^x \leq e^{A/e}$ for $A,x>0$, we have 
    \begin{equation}
        S(\lambda, \mu) (2C)^{-k} \leq s^k \pg \frac{ek}{\ell_I}\pd^{\ell_I} \leq s^k e^{k}.
    \end{equation}
    \item If $\frac{n}{s\sqrt{k}} \geq e$, then the function $\ell \in [1, k] \mapsto \pg \frac{n}{s\sqrt{k}} \frac{k}{\ell}\pd^{\ell}$ reaches its maximum at $\ell = k$, and therefore we have
    \begin{equation}
        S(\lambda, \mu) (2C)^{-k} \leq s^k \pg \frac{n}{s\sqrt{k}}\pd^{k} = \pg \frac{n}{\sqrt{k}} \pd^k.
    \end{equation}
    \end{enumerate}
We conclude that in all cases we have
    \begin{equation}
        S(\lambda, \mu) \leq (2eC)^k \pg \max \pg s, \frac{n}{\sqrt{k}} \pd \pd^k. \qedhere
    \end{equation}
\end{proof}

\begin{proof}[Proof of Theorem \ref{thm: bound on skew dimensions intro}]
Recall from \eqref{eq: NHLF forme avec somme excitée} that $\frac{d_{\lambda\backslash\mu}}{d_\lambda} = \frac{S(\lambda, \mu)}{n^{\downarrow k}}$, and that $n^{\downarrow k} \geq (n/e)^k$.
The result then follows immediately from Proposition \ref{prop: borne sommes excitées forme générale}.
\end{proof}

\subsection{Sharpness of the bound}\label{s: sharpness of the bound}

In this section we show that Theorem \ref{thm: bound on skew dimensions intro} is sharp, by giving lower bounds on the quantity $\frac{d_{\lambda\backslash \mu}}{d_\lambda}$ when both $\lambda$ and $\mu$ are rectangles. 

Let $n$ be a large integer, $1\leq k \leq n$, and $\Tilde{s} = \Tilde{s}(n) \geq \sqrt{n}$. We may assume that $h:= n/\Tilde{s}$, $m:=\sqrt{k}$, and $\ell:= k/h$ are integers. In the two examples below $\lambda = [\Tilde{s}^h]$. Note that $\Tilde{s}\asymp s(\lambda) =:s$.

\begin{example}
 Assume that $\sqrt{k} \geq h$ and set $\mu = [\ell^h]$. Then $\cE(\lambda, \mu) = \ag \mu \ad$, so by the Naruse hook length formula we have
 \begin{equation}
     \frac{d_{\lambda\backslash \mu}}{d_\lambda} = \frac{H(\lambda, \mu)}{n^{\downarrow k}} \geq \frac{(\Tilde{s}^{\downarrow \ell})^h}{n^k} \geq \frac{(\Tilde{s}/e)^{\ell h}}{n^k} = e^{-k} \pg \frac{\Tilde{s}}{n}\pd^k \geq (2e)^{-k} \pg \frac{s}{n}\pd^k.
 \end{equation}
\end{example}

\begin{example}
    Assume that $\sqrt{k} \leq h$ and set $\mu = [m^m]$. Rather than applying the Naruse hook length formula, we apply a symmetry argument: we have $d_{\lambda\backslash \mu}=d_{\lambda \backslash \Tilde{\mu}}$, where $\Tilde{\mu}$ is the top-right-most $m \times m$ square of $\lambda$. Let also $\gamma$ be the rectangle below $\Tilde{\mu}$ (i.e.\ the set of boxes in the first $h-m$ rows and last $m$ columns) and $\beta$ be the rectangle on the left on $\Tilde{\mu}$ (i.e.\ the set of boxes in the first $h-m$ columns and last $m$ rows). Set finally $\nu = \lambda\backslash \Tilde{\mu}$. We illustrate the notation in Figure \ref{fig: example with rectangles}.

\begin{figure}[!ht]
     \begin{center}
        \begin{ytableau}
            \none &&&&&&&&&& \\
            \none &&&&&&&&&& \\
            \none &&&&&&&&&& \\
            \none &&&&&&&&&& \\
            \none &*(black)&*(black)&*(black)&&&&&&& \\
             \none &*(black)&*(black)&*(black)&&&&&&& \\
            \none &*(black)&*(black)&*(black)&&&&&&& &\none\\
        \end{ytableau}
       \begin{ytableau}
            \none &*(lightgray)&*(lightgray)&*(lightgray)&*(lightgray)&*(lightgray)&*(lightgray)&*(lightgray)&*(black)&*(black)&*(black) \\
            \none &*(lightgray)&*(lightgray)&*(lightgray)&*(lightgray)&*(lightgray)&*(lightgray)&*(lightgray)&*(black)&*(black)&*(black) \\
            \none &*(lightgray)&*(lightgray)&*(lightgray)&*(lightgray)&*(lightgray)&*(lightgray)&*(lightgray)&*(black)&*(black)&*(black) \\
            \none &&&&&&&&*(darkgray)&*(darkgray)&*(darkgray)\\
            \none &&&&&&&&*(darkgray)&*(darkgray)&*(darkgray) \\
            \none &&&&&&&&*(darkgray)&*(darkgray)&*(darkgray) \\
            \none &&&&&&&&*(darkgray)&*(darkgray)&*(darkgray) &\none \\
        \end{ytableau}
    \end{center}
\caption{On both diagrams, $\lambda$ is the set of all boxes. On the left, the diagram $\mu$ is in black. On the right the diagrams $\Tilde{\mu}$, $\gamma$, and $\beta$ are respectively in black, dark gray, and light gray.}
\label{fig: example with rectangles}
\end{figure}

Then by the hook length formula, since most hook lengths compensate, we have
\begin{equation}
\begin{split}
     \frac{d_{\lambda\backslash \mu}}{d_\lambda} = \frac{d_{\nu}}{d_\lambda}  = \frac{(n-k)!}{n!} \frac{H(\lambda, \lambda)}{H(\nu, \nu)}
    =\frac{e^{O(k)}}{n^k} \frac{H(\lambda, \lambda)}{H(\nu, \nu)}= \frac{e^{O(k)}}{n^k} H(\mu, \mu) \frac{H(\lambda, \gamma)}{H(\nu,\gamma)}\frac{H(\lambda, \beta)}{H(\nu, \beta)}
\end{split}
\end{equation}
and therefore, using at the end that $h\ell = n$, we have
\begin{equation}
    \frac{d_{\lambda\backslash \mu}}{d_\lambda} \geq \frac{e^{O(k)}}{n^k} m^k \frac{h^k}{m^k} \frac{\ell^k}{m^k} = e^{O(k)} \frac{1}{m^k}  =  \pg \frac{e^{O(1)}}{\sqrt{k}} \pd^k.
\end{equation}
\end{example}

These two examples show that the bound in Theorem \ref{thm: bound on skew dimensions intro} is sharp (up to the value of the constant $C$) for any orders of magnitude for $k$ and $s(\lambda)$.
In particular this shows that when we decomposed diagrams into unions of line (Figure \ref{fig:stairs slicing}) and computed the contribution of each line doing thick hook decompositions (Figure \ref{fig: minimally excited diagram with respect to a thick hook decomposition}) this only led in general to an $e^{O(k)}$ error term!

\section{Bounds on characters}\label{s:bounds on characters}

The most canonical way to compute the characters of symmetric groups is via the Murnaghan--Nakayama rule. This is an explicit formula for characters (explained below in this section) which involves sums over ribbon tableaux and sign cancellations that can be very delicate to understand for large diagrams. The Murnaghan--Nakayama rule led to fruitful \textit{exponential} character bounds \cite{Roichman1996, MullerSchlagePuchta2007precutoff, LarsenShalev2008}, i.e.\ bounds uniform on $\lambda\in \widehat{\kS_n}$ and $\sigma$ is in $\kS_n$ (or a large subset of $\kS_n$) of the form 
\begin{equation}\label{eq: exponential character bound definition}
|\chi^\lambda(\sigma)|\leq d_\lambda^{-(a+o(1))\frac{\ln(n/f(\sigma))}{\ln n}},
\end{equation}
where $f(\sigma)$ is the number of fixed points of $\sigma$ (or 1 if $\sigma$ is fixed point free), and $a>0$ is a constant to optimize. These character bounds are related to \textit{mixing times} of random walks on $\kS_n$, and more specifically to generalizations of results of Diaconis and Shahashahani \cite{DiaconisShahshahani1981} on random transpositions. 
Recently, the analytic technique of \textit{hypercontractivity} (see for instance \cite{KeevashLifshitzLongMinzer2024HypercontractivityGlobalFunctions}) was successfully used to obtain exponential character bounds \cite{LifshitzMarmor2024charactershypercontractivity}. These bounds essentially show that \eqref{eq: exponential character bound definition} holds with $a=1$ for permutations with at most $n/(\ln n)^{O(1)}$ cycles, which improved on the Larsen--Shalev bounds for these permutations. It turned out that improving some of the technical estimates of \cite{LarsenShalev2008} led to stronger bounds \cite[Theorem 1.7, see also Remark 5.1]{TeyssierThevenin2025virtualdegreeswitten}, which are suited for permutations with $o(n)$ fixed points. For permutations with order $n$ fixed points, we obtained the sharpest possible exponential character bounds in \cite{Olesker-TaylorTeyssierThevenin2025sharpboundsandcutoff}, namely that \eqref{eq: exponential character bound definition} holds uniformly with $a = 1$ and no $o(1)$ error term.

\medskip

The character bounds considered in the present paper follow the line of work of \cite{Kerov1998InterlacingMeasures, Biane1998RepresentationsFreeProbability, Biane2001ApproximateFactorizationConcentration, RattanSniady2008, FeraySniady2011}. They have a different form and are not comparable to exponential character bounds.
For exponential character bounds the Young diagrams that are limiting for the bound are those with a very long first row. Such shapes impose strong constraints on the ribbon tableaux and it makes it \textit{more} possible to handle them using the Murnaghan--Nakayama rule. 
However the Murnaghan--Nakayama rule seemed to be little suited for diagrams with short rows and columns (and in particular for balanced diagrams), where there are less constraints on ribbon tableaux: Féray and Śniady \cite{FeraySniady2011} described it as being “cumbersome and hence intractable for computing characters corresponding to large Young diagrams”. It might therefore be surprising that the proof of Theorem \ref{thm: borne avec racine de sigma cas général} (via that of Theorem \ref{thm: character bound useful for fixed point free intro}) relies on the Murnaghan--Nakayama rule, and leads to precise character bounds in some regimes not previously covered in \cite{RattanSniady2008} or \cite{FeraySniady2011}.

\subsection{Ribbon tableaux and the Murnaghan--Nakayama rule}

We follow the terminology of \cite[Section 3.1]{LivreMeliot2017RepresentationTheoryofSymmetricGroups} for the Murnaghan--Nakayama rule and the related notions of ribbon tableaux and height.

\subsubsection{Skew partitions and ribbons}
If $\lambda$ and $\mu$ are Young diagrams such that $\mu \subset \lambda$, we recall that the \textit{skew diagram} $\lambda\backslash \mu$ is defined as the diagram containing the boxes that are in $\lambda$ and not in $\mu$. In particular, it is not necessarily connected. We illustrate this in Figure \ref{fig: def skew diagram}.

\begin{figure}[!ht]
    \begin{center}
\begin{ytableau}
    \none   & \\
    \none   &&&& \\
    \none   &&&&& \\
    \none   &&&&&&& 
\end{ytableau}
\begin{ytableau}
    \none  \\
    \none  \\
    \none[\backslash] \\
    \none
\end{ytableau}
\begin{ytableau}
    \none \\
    \none   && \\
    \none   &&& &\none&\none&\none[=] \\
    \none   &&&&& 
\end{ytableau}
\begin{ytableau}
    \none   & \\
    \none   &\none&\none&& \\
    \none   &\none&\none&\none&& \\
    \none   &\none&\none&\none&\none&\none&& &\none
\end{ytableau}
\end{center}   
\caption{The skew diagram $\lambda\backslash \mu$, where $\lambda = (7,5,4,1)$ and $\mu = (5,3,2)$.}
\label{fig: def skew diagram}
\end{figure}

A \textit{ribbon} is a skew diagram which is connected and does not contain a $[2,2]$ subdiagram. We give examples of diagrams which are ribbons in Figure \ref{fig: def ribbon positive examples} and which are not ribbons in Figure \ref{fig: def ribbon negative examples}.

\begin{figure}[!ht]
    \begin{center}
\begin{ytableau}
    \none \\
    \none & \\
    \none &\\
    \none   
\end{ytableau}
\begin{ytableau}
    \none \\
    \none \\
    \none   && \\
    \none   &\none& 
\end{ytableau}
\begin{ytableau}
    \none \\
    \none \\
    \none   & \\
    \none   &&&&& 
\end{ytableau}
\begin{ytableau}
    \none   & \\
    \none   && \\
    \none   &\none&& \\
    \none   &\none&\none&&&& &\none
\end{ytableau}
\end{center}
\caption{Some examples of ribbons.}
\label{fig: def ribbon positive examples}
\end{figure}

\begin{figure}[!ht]
    \begin{center}
\begin{ytableau}
    \none \\
    \none &&*(red)\\
    \none   &\none&\none&*(red)& \\
    \none   &\none&\none&\none&& 
\end{ytableau}
\begin{ytableau}
    \none   & \\
    \none   && \\
    \none   &\none&*(red)&*(red) \\
    \none   &\none&*(red)&*(red)&&&
\end{ytableau}
\begin{ytableau}
    \none   &&&&\\
    \none   &&\none&\none&&\\
    \none   &&\none&\none&\none& \\
    \none   &&\none&\none&\none&&& &\none
\end{ytableau}
\end{center}
\caption{Some examples of skew diagrams that are not ribbons. The first one is disconnected (the part where it is disconnected is emphasized in red), the second one contains a $[2,2]$ subdiagram (emphasized in red), and the third one is not a skew diagram.}
\label{fig: def ribbon negative examples}
\end{figure}

A standard tableau of a diagram (a Young diagram or a skew diagram) of size $n$ is a way to place the numbers from $1$ to $n$ in the boxes of the diagram, in such a way that numbers increase along each row and each column. We denote the set of standard tableaux of a diagram $\lambda$ by $\ST(\lambda)$, and by $\ST(\lambda\backslash \mu)$ that of a skew diagram $\lambda\backslash \mu$.

\subsubsection{Ribbon tableaux and height}
A \textit{ribbon tableau} of a partition $\lambda$ is a sequence of partitions $\emptyset = \mu^{(0)} \subset \mu^{(1)} \subset \ldots \subset \mu^{(N)} = \lambda$ such that for each $i \in \{1,\ldots, N\}$, $r^{(i)}:= \mu^{(i)}\backslash \mu^{(i-1)}$ is a ribbon. Note that it is important to remember not only the shape of the ribbons $r^{(i)}$'s, but also their position in $\lambda$. We say that a ribbon tableau of $\lambda$ has \textit{weight} $\alpha = (\alpha_1, \alpha_2,\ldots, \alpha_N)$ if $\bg r^{(i)} \bd = \alpha_i$ for $i\geq 1$, where $|E|$ denotes the number of boxes in a set of boxes $E$.
We denote by $\RT(\lambda, \alpha)$ the set of ribbon tableaux of shape $\lambda$ and weight $\alpha$, and $\RT(\lambda):= \bigcup_\alpha \RT(\lambda,\alpha)$ the set of all ribbon tableaux of shape $\lambda$. We illustrate this in Figure \ref{fig: def ribbon tableau}.

\begin{figure}[!ht]
    \begin{center}

\begin{ytableau}
    \none   &\\
    \none   &&&&\\
    \none   &&&&&\\
    \none   &&&&&&&
\end{ytableau}
\begin{ytableau}
    \none   &\\
    \none   &&&&\\
    \none   &&&&&*(black)\\
    \none   &&&&&&&
\end{ytableau}
\begin{ytableau}
    \none   &\\
    \none   &&&&*(black)\\
    \none   &&&&\\
    \none   &&&&&&&
\end{ytableau}

\begin{ytableau}
    \none &\none\\
    \none   &\\
    \none   &&*(black)&*(black)\\
    \none   &&&&\\
    \none   &&&&&&&
\end{ytableau}
\begin{ytableau}
    \none &\none\\
    \none   &\\
    \none   &\\
    \none   &&&&\\
    \none   &&&&&&*(black)&*(black)
\end{ytableau}
\begin{ytableau}
    \none &\none\\
    \none   &\\
    \none   &\\
    \none   &&&*(black)&*(black)\\
    \none   &&&&&&\none&\none
\end{ytableau}

\begin{ytableau}
    \none &\none\\
    \none   &*(black)\\
    \none   &*(black)\\
    \none   &*(black)&*(black)\\
    \none   &&&&&&\none&\none
\end{ytableau}
\begin{ytableau}
    \none &\none\\
    \none &\none\\
    \none &\none\\
    \none &\none\\
    \none   &*(black)&*(black)&*(black)&*(black)&*(black)&\none&\none
\end{ytableau}
\begin{ytableau}
    \none &\none\\
    \none &\none\\
    \none &\none\\
    \none &\none &\none &\none[\emptyset]\\
    \none   &\none&\none&\none&\none&\none&\none&\none
\end{ytableau}
\end{center}
\caption{A representation of the ribbon tableau corresponding to the sequence $\emptyset \subset [5]  \subset [5,2,1,1] \subset [5,4,1,1] \subset [7,4,1,1] \subset  [7,4,3,1] \subset [7,4,4,1] \subset [7,5,4,1],$ which is a ribbon tableau of $\lambda = [7,5,4,1]$ with weight $\alpha = [5,4,2,2,2,1,1]$. This is represented by starting from $\lambda$ and peeling the ribbons (represented in black) one by one.}
\label{fig: def ribbon tableau}
\end{figure}

Finally, we can represent ribbon tableaux in the diagram itself, writing the integer $i$ in each ribbon $r^{(i)}$. This is a more compact and more visual representation, which we use from now on. We illustrate this in Figure \ref{fig: ribbon tableau visual representation}. 

\begin{figure}[!ht]
    \begin{center}
\begin{ytableau}
    \none   &*(orange!60)2\\
    \none   &*(orange!60)2&*(cyan!60)5&*(cyan!60)5&*(blue!60)6\\
    \none   &*(orange!60)2&*(orange!60)2&*(yellow!60)3&*(yellow!60)3&*(violet!60)7\\
    \none   &*(red!60)1&*(red!60)1&*(red!60)1&*(red!60)1&*(red!60)1&*(green!60)4&*(green!60)4&\none
\end{ytableau}
    \begin{ytableau}
    \none   &*(orange!60)2\\
    \none   &*(orange!60)2&*(cyan!60)5&*(cyan!60)5&*(violet!60)7\\
    \none   &*(orange!60)2&*(orange!60)2&*(green!60)4&*(green!60)4&*(blue!60)6\\
    \none   &*(red!60)1&*(red!60)1&*(red!60)1&*(red!60)1&*(red!60)1&*(yellow!60)3&*(yellow!60)3&\none
\end{ytableau}
 \begin{ytableau}
    \none   &*(blue!60)6\\
    \none   &*(red!60)1&*(green!60)4&*(green!60)4&*(violet!60)7\\
    \none   &*(red!60)1&*(orange!60)2&*(orange!60)2&*(orange!60)2&*(yellow!60)3\\
    \none   &*(red!60)1&*(red!60)1&*(red!60)1&*(orange!60)2&*(yellow!60)3&*(cyan!60)5&*(cyan!60)5&\none
\end{ytableau}
\end{center}
\caption{Three ribbon tableaux of $\lambda = [7,5,4,1]$ with weight $\alpha = [5,4,2,2,2,1,1]$, corresponding, from left to right, to the sequences $\emptyset \subset [5] \subset [5,2,1,1] \subset [5,4,1,1] \subset [7,4,1,1] \subset  [7,4,3,1] \subset [7,4,4,1] \subset [7,5,4,1]$ (the same as in Figure \ref{fig: def ribbon tableau}), to $\emptyset \subset [5]  \subset [5,2,1,1] \subset [7,2,1,1] \subset [7,4,1,1] \subset  [7,4,3,1] \subset [7,5,3,1] \subset [7,5,4,1]$, and to $\emptyset \subset [3,1,1]  \subset [4,4,1] \subset [5,5,1] \subset [5,5,3] \subset  [7,5,3] \subset [7,5,3,1] \subset [7,5,4,1]$. Importantly, the ribbons of the first two ribbon tableaux are the same, but since the order in which they are peeled matters, they lead to two distinct ribbon tableaux. 
Also in this diagram and all the following ones, the colors play the same role as numbers, and are therefore redundant.}
\label{fig: ribbon tableau visual representation}
\end{figure}

The difference between the vertical coordinates of the highest box and lowest box of a ribbon $r$ is called \textit{height} (of $r$), and is denoted by $\height(r)$. See Figure \ref{fig: def height of ribbons} for examples.

\begin{figure}[!ht]
       
\begin{center}
\begin{ytableau}
    \none \\
    \none  \\
    \none \\
    \none &&&  
\end{ytableau}
\begin{ytableau}
    \none \\
    \none \\
    \none   && \\
    \none   &\none& 
\end{ytableau}
\begin{ytableau}
    \none \\
    \none &\\
    \none   & \\
    \none   &&&&& 
\end{ytableau}
\begin{ytableau}
    \none   & \\
    \none   && \\
    \none   &\none&& \\
    \none   &\none&\none&&&& &\none[.]
\end{ytableau}
\end{center}
\caption{Examples of ribbons with respective heights (from left to right) 0, 1, 2, and~3.}
\label{fig: def height of ribbons}
\end{figure}

If $\lambda$ is a partition and $T=(r^{(1)}, \ldots, r^{(N)})\in \RT(\lambda)$ is a ribbon tableau of $\lambda$, we define the \textit{height} of $T$ as the sum of the heights of all ribbons in $T$: $\height(T) = \sum_{i=1}^N \height(r^{(i)})$. See Figure \ref{fig: def height of ribbon tableaux} for examples.

\begin{figure}[!ht]
    \begin{center}
        \begin{ytableau}
    \none \\
    \none  \\
    \none \\
    \none &*(red!60)1&*(red!60)1&*(red!60)1  
\end{ytableau}
\begin{ytableau}
    \none \\
    \none  \\
    \none &*(yellow!60)3\\
    \none &*(red!60)1&*(red!60)1&*(orange!60)2  
\end{ytableau}
\begin{ytableau}
    \none \\
    \none  \\
    \none &*(red!60)1\\
    \none &*(red!60)1&*(red!60)1&*(orange!60)2  
\end{ytableau}
\begin{ytableau}
    \none \\
    \none &*(red!60)1\\
    \none &*(red!60)1&*(orange!60)2\\
    \none &*(red!60)1&*(orange!60)2&*(orange!60)2
\end{ytableau}
\begin{ytableau}
    \none &*(red!60)1 \\
    \none &*(red!60)1&*(cyan!60)5\\
    \none &*(red!60)1&*(orange!60)2&*(green!60)4&*(blue!60)6\\
    \none &*(red!60)1&*(orange!60)2&*(yellow!60)3&*(yellow!60)3  
\end{ytableau}
    \end{center}
\caption{Examples of ribbon tableaux, with respective heights $0$, $0+0+0=0$, $1+0=1$, $2+1=3$, and $3+1+0+0+0+0=4$.}
\label{fig: def height of ribbon tableaux}
\end{figure}

\subsubsection{The Murnaghan--Nakayama rule}\label{s: MN rule}

The Murnaghan--Nakayama rule (\cite[Theorem 3.10]{LivreMeliot2017RepresentationTheoryofSymmetricGroups}) is a fascinating formula, which not only shows that the characters $\ch^\lambda(\sigma)$ are actually integers, but also provides us with a way to compute them. It can be stated as follows.

\begin{theorem}[Murnaghan-Nakayama rule]
\label{theorem:MN}
Let $n\geq 1$, $\lambda\vdash n$ be a partition of $n$, and $\sigma\in \mathfrak{S}_n$. 
Denote by $\alpha$ the cycle type of $\sigma$, that is, the list of cycle lengths of $\sigma$ sorted in decreasing order. Then, we have
\begin{equation}\label{e:MN}
    \ch^\lambda(\sigma) = \sum_{T \in \RT(\lambda, \alpha)} (-1)^{\height(T)}. 
\end{equation}
\end{theorem}

Since the character $\ch^\lambda(\sigma)$ depends only on the cycle type $\alpha$ of $\sigma$, we may write $\ch^\lambda(\alpha)$ instead of $\ch^\lambda(\sigma)$.

A useful property of the Murnaghan--Nakayama rule is that the value of the sum does not depend on the order of the entries of $\alpha$ (see \cite[Corollary 7.17.4 and the comment below the proof of Proposition 7.17.6]{LivreStanleyVol2}, or the proof of \cite[Proposition 3.10]{LivreMeliot2017RepresentationTheoryofSymmetricGroups}).
When computing characters, it can be convenient to permute some entries of $\alpha$ to reduce the number of diagrams to consider. Let us illustrate this on two examples.

\begin{example}\label{ex:MN 32 311 et 131 et 113}
Let $\lambda = [3,2]$ and $\sigma$ be the 3-cycle $(1\, 4 \, 5)$ in $\mathfrak{S}_5$. The cycle type of $\sigma$ is then $\alpha =(3,1,1)$, but can be also reordered as $\alpha = (1,3,1)$ and $\alpha = (1,1,3)$.

    The set $\RT(\lambda, (3,1,1))$ consists of 3 ribbon tableaux (see Figure \ref{fig: RT 32 311}), while $\RT(\lambda, (3,1,1))$ both consist of a single tableau (see Figure \ref{fig: RT 32 131 and 113}). Applying the Murnaghan--Nakayama rule \eqref{e:MN} with $\alpha = (3,1,1)$, we obtain $  \ch^\lambda(\sigma) = (-1)^0 + (-1)^1 + (-1)^1 = -1$, while if we take $\alpha \in \ag (1,3,1), (1,1,3) \ad$ we immediately obtain $\ch^\lambda(\sigma) = (-1)^1 = -1$.

\begin{figure}[!ht]
        
\begin{center}
\begin{ytableau}
    \none &*(orange!60)2&*(yellow!60)3\\
    \none &*(red!60)1&*(red!60)1&*(red!60)1  &\none&\none
\end{ytableau}
\begin{ytableau}
    \none &*(red!60)1&*(orange!60)2\\
    \none &*(red!60)1&*(red!60)1  &*(yellow!60)3 &\none&\none
\end{ytableau}
\begin{ytableau}
    \none &*(red!60)1&*(yellow!60)3\\
    \none &*(red!60)1&*(red!60)1  &*(orange!60)2&\none
\end{ytableau}
\end{center}
\caption{The three elements of $\RT(\lambda, (3,1,1))$, where $\lambda = [3,2]$. They have respective heights (from left to right) $0+0+0=0$, $1+0+0=1$, and $1+0+0=1$.}
\label{fig: RT 32 311}
\end{figure}

\begin{figure}[!ht]
        
\begin{center}
\begin{ytableau}
    *(orange!60)2&*(orange!60)2\\
    *(red!60)1&*(orange!60)2&*(yellow!60)3 &\none &\none &\none &\none
\end{ytableau}
\begin{ytableau}
    *(orange!60)2&*(yellow!60)3\\
    *(red!60)1&*(yellow!60)3&*(yellow!60)3
\end{ytableau}
\end{center}
\caption{Here $\lambda = [3,2]$. On the left: the single element of $\RT(\lambda, (1,3,1))$, which has height $0+1+0 = 1$. On the right: the single element of $\RT(\lambda, (1,1,3))$, which has height $0+0+1 = 1$. }
\label{fig: RT 32 131 and 113}
\end{figure}

\end{example}


\begin{example}\label{ex:MN 433 33211}
Let $\lambda = [4,3,3]$, and let $\sigma$ with cycle type $\alpha = (3,3,2,1,1)$. If we compute $\ch^\lambda(\sigma)$ without reordering $\alpha = (3,3,2,1,1)$, we have to consider the 12 ribbon tableaux in $\RT(\lambda, \alpha)$ (see Figure \ref{fig: twelve ribbon tableaux}). 7 of these tableaux have an even sign (and 5 an odd sign) which leads to $\ch^\lambda(\sigma) = 7-5 = 2$. 
On the other hand, if we reorder $\alpha = (3,3,2,1,1)$ as $\beta = (1,3,3,2,1)$, we only have to consider the 2 ribbon tableaux in $\RT(\lambda, \beta)$ (see Figure \ref{fig: two ribbon tableaux}), which both have an even sign, and we immediately obtain $\ch^\lambda(\sigma) = 2$.

\begin{figure}[!ht]
    
\begin{center}
\begin{ytableau}
    \none &*(yellow!60)3&*(yellow!60)3&*(cyan!60)5\\
    \none &*(orange!60)2&*(orange!60)2&*(orange!60)2\\
    \none &*(red!60)1&*(red!60)1&*(red!60)1&*(green!60)4
\end{ytableau}
\begin{ytableau}
    \none &*(yellow!60)3&*(yellow!60)3&*(green!60)4\\
    \none &*(orange!60)2&*(orange!60)2&*(orange!60)2\\
    \none &*(red!60)1&*(red!60)1&*(red!60)1&*(cyan!60)5
\end{ytableau}
\begin{ytableau}
    \none &*(yellow!60)3&*(yellow!60)3&*(cyan!60)5\\
    \none &*(red!60)1&*(orange!60)2&*(orange!60)2\\
    \none &*(red!60)1&*(red!60)1&*(orange!60)2&*(green!60)4
\end{ytableau}
\begin{ytableau}
    \none &*(yellow!60)3&*(yellow!60)3&*(green!60)4\\
    \none &*(red!60)1&*(orange!60)2&*(orange!60)2\\
    \none &*(red!60)1&*(red!60)1&*(orange!60)2&*(cyan!60)5
\end{ytableau}

\begin{ytableau}
    \none \\
    \none &*(orange!60)2&*(orange!60)2&*(cyan!60)5\\
    \none &*(red!60)1&*(orange!60)2&*(green!60)4\\
    \none &*(red!60)1&*(red!60)1&*(yellow!60)3&*(yellow!60)3
\end{ytableau}
\begin{ytableau}
    \none \\
    \none &*(orange!60)2&*(orange!60)2&*(cyan!60)5\\
    \none &*(red!60)1&*(orange!60)2&*(yellow!60)3\\
    \none &*(red!60)1&*(red!60)1&*(yellow!60)3&*(green!60)4
\end{ytableau}
\begin{ytableau}
    \none \\
    \none &*(orange!60)2&*(orange!60)2&*(green!60)4\\
    \none &*(red!60)1&*(orange!60)2&*(yellow!60)3\\
    \none &*(red!60)1&*(red!60)1&*(yellow!60)3&*(cyan!60)5
\end{ytableau}
\begin{ytableau}
    \none \\
    \none &*(red!60)1&*(green!60)4&*(cyan!60)5\\
    \none &*(red!60)1&*(yellow!60)3&*(yellow!60)3\\
    \none &*(red!60)1&*(orange!60)2&*(orange!60)2&*(orange!60)2
\end{ytableau}

\begin{ytableau}
    \none \\
    \none \\
    \none \\
    \none
\end{ytableau}
\begin{ytableau}
    \none \\
    \none &*(red!60)1&*(yellow!60)3&*(cyan!60)5\\
    \none &*(red!60)1&*(yellow!60)3&*(green!60)4\\
    \none &*(red!60)1&*(orange!60)2&*(orange!60)2&*(orange!60)2
\end{ytableau}
\begin{ytableau}
    \none \\
    \none &*(red!60)1&*(orange!60)2&*(cyan!60)5\\
    \none &*(red!60)1&*(orange!60)2&*(green!60)4\\
    \none &*(red!60)1&*(orange!60)2&*(yellow!60)3&*(yellow!60)3
\end{ytableau}
\begin{ytableau}
    \none \\
     \none &*(red!60)1&*(orange!60)2&*(cyan!60)5\\
    \none &*(red!60)1&*(orange!60)2&*(yellow!60)3\\
    \none &*(red!60)1&*(orange!60)2&*(yellow!60)3&*(green!60)4
\end{ytableau}
\begin{ytableau}
    \none \\
     \none &*(red!60)1&*(orange!60)2&*(green!60)4\\
    \none &*(red!60)1&*(orange!60)2&*(yellow!60)3\\
    \none &*(red!60)1&*(orange!60)2&*(yellow!60)3&*(cyan!60)5
\end{ytableau}
\begin{ytableau}
    \none \\
    \none \\
    \none \\
    \none
\end{ytableau}
\end{center}
\caption{The twelve ribbon tableaux in $\RT(\lambda, \alpha)$, where $\lambda = [4,3,3]$ and $\alpha = (3,3,2,1,1)$. They have respective heights 0, 0, 2, 2, 2, 3, 3, 2, 3, 4, 5, and 5.}
\label{fig: twelve ribbon tableaux}
\end{figure}

\begin{figure}[!ht]
    
\begin{center}
\begin{ytableau}
    \none &*(green!60)4&*(green!60)4&*(cyan!60)5\\
    \none &*(yellow!60)3&*(yellow!60)3&*(yellow!60)3\\
    \none &*(red!60)1&*(orange!60)2&*(orange!60)2&*(orange!60)2&\none&\none
\end{ytableau}
\begin{ytableau}
    \none &*(green!60)4&*(green!60)4&*(cyan!60)5\\
    \none &*(orange!60)2&*(orange!60)2&*(yellow!60)3\\
    \none &*(red!60)1&*(orange!60)2&*(yellow!60)3&*(yellow!60)3
\end{ytableau}
\end{center}
\caption{The two ribbon tableaux in $\RT(\lambda, \beta)$, where $\lambda = [4,3,3]$ and $\beta = (13,3,2,1)$. They have respective heights 0 and 2, which are both even.}
\label{fig: two ribbon tableaux}
\end{figure}
\end{example}

The Murnaghan--Nakayama rule is actually fairly simple to use for small diagrams, but using it to estimate characters of large diagrams can be very complicated.

There are two main difficulties to overcome. The first one is that there could be significant cancellations between ribbon tableaux since the summands $(-1)^{\height(T)}$ can be positive or negative. A possibility to handle this is to forget about compensations by applying the triangle inequality in \eqref{e:MN}, which leads to the bound

\begin{equation} \label{eq: bound characters triangle inequality number of ribbon tableaux}
     \bg \ch^\lambda(\sigma) \bd = \bg \sum_{T \in \RT(\lambda, \alpha)} (-1)^{\height(T)} \bd \leq \bg \RT(\lambda, \alpha) \bd.
\end{equation}
One then only has to estimate the number of ribbon tableaux of $\lambda$ of weight $\alpha$. Applying the triangle inequality in the Murnaghan--Nakayama rule \eqref{e:MN} seems wasteful at first glance but can be surprisingly precise in practice.
The second and main difficulty is that when permutations have many small cycles and diagrams do not have a long first row, the number of associated ribbon tableaux is large and difficult to control. 

\subsubsection{A character identity seen via the Murnaghan--Nakayama rule}

There are many identities involving characters. One of them is the iterated branching rule written below as Lemma \ref{lem: MN rewritten with sigma star}, which was used for instance in the proof of Lemma 10 in \cite{DousseFeray2019skewdiagramscharacters}. We observed in \cite{Olesker-TaylorTeyssierThevenin2025sharpboundsandcutoff} that this identity is useful to leverage character bounds for fixed point free permutations to permutations of any support size, and will also use this to prove Theorem \ref{thm: borne avec racine de sigma cas général}. 

Let $n$ and $k$ be integers such that $2\leq k\leq n$, and let $\sigma$ is a permutation of $\kS_n$ with support $k$. Let $\sigma^* \in \kS_k$ be a permutation with as many cycles of length $i$ as $\sigma$ for $i\geq 2$, but no fixed point; to define it uniquely we fill the cycles of $\sigma^*$ in (weakly) decreasing order of their lengths with the numbers from 1 to $k$.
For example if $\sigma = ( 1\; 3 \; 5 \; 8) (2 \; 7)(4)(6)(9) \in \kS_9$, then $\sigma^* = (1 \; 2 \; 3 \; 4)(5 \; 6) \in \kS_6$.

A way to view the following iterated branching rule is that it corresponds to peeling first the $n-k$ ribbons of size 1 in the Murnaghan--Nakayama rule.

\begin{lemma}\label{lem: MN rewritten with sigma star}
    Let $2\leq k \leq n$, $\sigma\in \kS_n$ with support of size $k$, and $\lambda \vdash n$. Then
    \begin{equation}
        \ch^\lambda(\sigma) = \sum_{\mu \subset_{\vdash k} \lambda} \ch^\mu(\sigma^*) d_{\lambda\backslash \mu}.
    \end{equation}
Equivalently, we have
    \begin{equation}
        \chi^\lambda(\sigma) = \sum_{\mu \subset_{\vdash k} \lambda} \ch^\mu(\sigma^*) \frac{d_{\lambda\backslash \mu}}{d_\lambda}.
    \end{equation}
\end{lemma}
\begin{proof}
Let $\ell$ be the length of the longest cycle in $\sigma$. We may assume without loss of generality that $\ell \geq 2$. Let $\alpha_i$ be the number of $i$-cycles in $\sigma$, and set $\alpha = \pg 1^{\alpha_1}, \ell^{\alpha_\ell}, (\ell-1)^{\alpha_{\ell-1}}, \ldots, 2^{\alpha_2}\pd $, the cycle type of $\sigma$, written with fixed points first and the other cycles in (weakly) decreasing order of their lengths. Let also $\alpha^* = \pg \ell^{\alpha_\ell}, (\ell-1)^{\alpha_{\ell-1}}, \ldots, 2^{\alpha_2}\pd $ be the cycle type of $\sigma^*$, where the cycles are written in (weakly) decreasing order of their lengths. Then by the Murnaghan--Nakayama rule, peeling the $\alpha_1$ ribbons of length 1, we have
\begin{equation}
    \ch^\lambda(\sigma) = \sum_{T \in \RT(\lambda, \alpha)} (-1)^{\height(T)} = \sum_{\mu \subset_{\vdash k} \lambda} \pg \sum_{T'\in \RT(\mu, \alpha^*)}(-1)^{\height(T')}\pd d_{\lambda\backslash \mu} = \sum_{\mu \subset_{\vdash k} \lambda} \ch^\mu(\sigma^*) d_{\lambda\backslash \mu}.
\end{equation}
The second part of the lemma is obtained from the first one by dividing both sides by $d_\lambda$.
\end{proof}

\subsection{Elementary character bounds in terms of the number of cycles}\label{s: Character bound in function of the diagonal length and the number of cycles}

Larsen and Shalev \cite{LarsenShalev2008} developed an elaborate inductive peeling process to prove exponential character bounds which are especially precise for fixed point free permutations and representations with a long first row. In this section show that applying a simpler peeling process that relies only on counting arguments leads to more precise bounds for diagrams that do not have a long first row, and in particular for balanced diagrams. 

\medskip

If $\sigma\in \kS_n$ and $j\geq 1$, $\cyc_j(\sigma)$ is the number of cycles of length $j$ of $\sigma$, $\cyc_{\geq j}(\sigma) = \sum_{i\geq j} \cyc_i(\sigma)$ is the number of cycles of length at least $j$ of $\sigma$, and $\cyc(\sigma) = \cyc_{\geq 1}(\sigma)$ is the total number of cycles of $\sigma$.
Then $\sum_{j\geq 1} j\cyc_j(\sigma) = n$ and $\sum_{j\geq 2} j\cyc_j(\sigma) = \supp(\sigma)$ is the support size of $\sigma$.
Also, for $\lambda \vdash n$, $\delta(\lambda)$ is the diagonal length of $\lambda$, that is, $\delta(\lambda) = \max \ag i\geq 1 \mid \min(\lambda_i, \lambda_i')\geq i \ad$, and $\ell(\lambda)$ is the number of rows of $\lambda$, that is, $\ell(\lambda) = \lambda_1'$.



\begin{lemma}\label{lem: elementary bound number of peelable ribbons}
Let $\lambda$ be a Young diagram. Let $r\geq 1$. 
\begin{enumerate}
    \item There are at most $\ell(\lambda)$ ways to peel a ribbon of length $r$ from $\lambda$.
    \item There are at most $2\delta(\lambda)$ ways to peel a ribbon of length $r$ from $\lambda$.
\end{enumerate}
\end{lemma}
\begin{proof}
(a) A ribbon, given its length, is determined by the position of its bottom-right-most box, which must be at the end of a row of $\lambda$. This is illustrated in Figure \ref{fig: young diagram brm possible positions of ribbons}.
Hence the number of ribbons of length $r$ of $\lambda$ is at most its number of rows, which is $\ell(\lambda)$.

\begin{figure}[!ht]
    \centering
{
\ytableausetup{boxsize=1.2em}
\begin{ytableau}
    \none &*(darkgray) \\
    \none &*(darkgray) \\
    \none &&&*(darkgray) \\
    \none &&&*(darkgray) \\
    \none &&&*(darkgray)  \\
    \none &&&&&*(darkgray)  \\
    \none &&&&&&&&*(darkgray) \\
    \none &&&&&&&&*(darkgray) \\
    \none &&&&&&&&&&&&*(darkgray) \\
    \none &&&&&&&&&&&&&&&*(darkgray) \\
    \none &&&&&&&&&&&&&&&&&&&&&&&&*(darkgray) \\
\end{ytableau}
}
    \caption{The Young diagram $\lambda = [24,15,12,8,8,5,3,3,3,1,1]$. The right-most box of each row is in dark gray. These boxes are the bottom-right-most possible positions of boxes of ribbons of $\lambda$.}
    \label{fig: young diagram brm possible positions of ribbons}
\end{figure}

(b) Pushing further the previous argument, a ribbon must satisfy at least one of the following two conditions: having its bottom-right-most box at the end of one of the $\delta(\lambda)$ first rows of $\lambda$; or having its top-left-box at the end of one of the first $\delta(\lambda)$ columns of $\lambda$. These boxes are illustrated in Figure \ref{fig: young diagram brm and tlm possible positions of ribbons}.
\begin{figure}[!ht]
    \centering
{
\ytableausetup{boxsize=1.2em}
\begin{ytableau}
    \none &*(lightgray) \\
    \none &\\
    \none &&*(lightgray)&*(lightgray)\\
    \none &&& \\
    \none &&&  \\
    \none &&&&*(lightgray)&*(lightgray)  \\
    \none &&&&&&&&*(darkgray) \\
    \none &&&&&&&&*(darkgray) \\
    \none &&&&&&&&&&&&*(darkgray) \\
    \none &&&&&&&&&&&&&&&*(darkgray) \\
    \none &&&&&&&&&&&&&&&&&&&&&&&&*(darkgray) \\
\end{ytableau}
}
    \caption{The Young diagram $\lambda = [24,15,12,8,8,5,3,3,3,1,1]$. Here $\delta(\lambda) = 5$. The right-most box of each of the first $\delta(\lambda)$ rows is in dark gray. The top box of each of the first $\delta(\lambda)$ columns is in light gray.}
    \label{fig: young diagram brm and tlm possible positions of ribbons}
\end{figure}
\end{proof}

We can now prove the desired character bound.

\begin{proof}[Proof of Theorem \ref{thm: character bound useful for fixed point free intro}]
   Let $\alpha$ be any ordering of the cycle lengths of $\sigma$, for example in (weakly) increasing order. Then applying the triangle inequality in the Murnaghan--Nakayama rule gives $|\ch^\lambda(\sigma)| \leq |\RT(\lambda, \alpha)|$, as written in \eqref{eq: bound characters triangle inequality number of ribbon tableaux}.

    We therefore only need to bound $|\RT(\lambda, \alpha)|$. We note that peeling a ribbon $\rho$ from $\lambda$ cannot increase its number of rows or its diagonal length. Hence, applying Lemma \ref{lem: elementary bound number of peelable ribbons} iteratively to the $\cyc(\sigma)$ ribbons that we peel gives $|\RT(\lambda, \alpha)|\leq \ell(\lambda)^{\cyc(\sigma)}$ and $|\RT(\lambda, \alpha)|\leq (2\delta(\lambda))^{\cyc(\sigma)}$, which concludes the proof.
\end{proof}

\subsection{Bounds on characters for general shapes of diagrams}\label{s: proof main character bound}

We are now ready to prove Theorem \ref{thm: borne avec racine de sigma cas général}.

\begin{proof}[Proof of Theorem \ref{thm: borne avec racine de sigma cas général}]
    Let $n$ be a large integer. Let $\lambda \vdash n$ and $\sigma\in \kS_n\backslash\ag \Id \ad$. Denote the support size of $\sigma$ by $k$. By Lemma \ref{lem: MN rewritten with sigma star}, the triangle inequality, and since the number of partitions of the integer $k$ is $e^{O(k)}$ (actually even $e^{O(\sqrt{k})}$), we have 
    \begin{equation}
        \bg \chi^\lambda(\sigma) \bd \leq \sum_{\mu\subset \lambda \du |\mu| = k} \bg \ch^\mu(\sigma^*)\bd \frac{d_{\lambda\backslash \mu}}{d_\lambda} \leq e^{O(k)} \max_{\mu\subset \lambda \du |\mu| = k} \bg \ch^\mu(\sigma^*)\bd \frac{d_{\lambda\backslash \mu}}{d_\lambda}.
    \end{equation}
We deduce from Theorem \ref{thm: character bound useful for fixed point free intro} (using that for any $\mu \vdash k$ the diagonal length of $\mu$ satisfies $\delta(\mu)\leq \sqrt{k}$) and Theorem \ref{thm: bound on skew dimensions intro} that
\begin{equation}
    \bg \chi^\lambda(\sigma) \bd \leq e^{O(k)} \sqrt{k}^{\cyc(\sigma^*)} \pg \max \pg  \frac{1}{\sqrt{k}} , \frac{s}{n}\pd \pd^k = e^{O(k)} \pg\frac{1}{\sqrt{k}}\pd^{k-\cyc(\sigma^*)} \max \pg  1 , \frac{s\sqrt{k}}{n}\pd^k.
\end{equation}
This concludes the proof since $|\sigma| \asymp k$, $e^{O(k)} = O(1)^k$, and $k-\cyc(\sigma^*) = |\sigma^*| = |\sigma|$.
\end{proof}

\section{Some comments on random compression probabilities}

Let $n$ and $k$ be two integers such that $1\leq k \leq n$. Let $\lambda\vdash n$ and $\mu\vdash k$ be two diagrams such that $\mu \subset \lambda$. 
A natural quantity to consider is the probability  
\begin{equation}
    P_\lambda(\mu) = \frac{d_\mu d_{\lambda \backslash \mu}}{d_\lambda}
\end{equation}
that the first $k$ integers fill $\mu$ in a random tableau of $\lambda$. This can also be interpreted as a random compression probability, see the survey \cite{Meliot2019RandomCompressions}. Denote the Plancherel measure evaluated at a Young diagram $\nu\vdash k$ by 
\begin{equation}
    \Pl_k(\nu) = \frac{d_\nu^2}{k!}.
\end{equation}
Dousse and Féray \cite{DousseFeray2019skewdiagramscharacters} studied the asymptotics of the quantity
\begin{equation}
    A_{\lambda, \mu} := k!\frac{d_{\lambda\backslash\mu}}{d_\mu d_\lambda} = \frac{P_\lambda(\mu)}{\Pl_k(\mu)},
\end{equation}
when $\lambda$ and $\mu$ are balanced Young diagrams, and proved precise bounds when $k=O(\sqrt{n})$.

An interesting question is to understand how \textit{close} the two measures $P_\lambda(\cdot)$ and $\Pl_k(\cdot)$ are, for different shapes of $\lambda$ and sizes $k$, and different notions of distances.
The quantity $\sum_{\mu\vdash k}|P_\lambda(\mu) - \Pl_k(\mu)|$ is (twice) the total variation distance between the two measures, while
\begin{equation}
    \max_{\mu \vdash k} |A_{\lambda, \mu} - 1|
\end{equation}
is a form of $L^\infty$ distance between $P_\lambda(\cdot)$ and $\Pl_k(\cdot)$. (If we view the Plancherel measure as the stationary measure of a Markov chain, this is how the $L^\infty$ distance from stationarity defined, see \cite[Section 4.7]{LivreLevinPeres2019MarkovChainsAndMixingTimesSecondEdition}.) Dousse and Féray restricted their attention to diagrams $\lambda$ and $\mu$ that are balanced. Their results are therefore bounds of $L^\infty$-type restricted to balanced diagrams.

\medskip

Dousse and Féray \cite[Conjecture 6]{DousseFeray2019skewdiagramscharacters} conjectured that for any $C>0$ there exists $C'>0$ such that as $n\to \infty$, uniformly over all $1\leq k\leq n$, and all $C$-balanced diagrams $\lambda\vdash n$ and $\mu\vdash k$ such that $\mu\subset \lambda$, we have
\begin{equation}
\label{eq: conj Dousse Feray}
    -C'\frac{k^{3/2}}{n^{1/2}} \leq \ln A_{\lambda, \mu} \leq C'\frac{k^{3/2}}{n^{1/2}}. 
\end{equation}
They proved \eqref{eq: conj Dousse Feray} for $k= o(n^{1/3})$ in \cite[Corollary 2]{DousseFeray2019skewdiagramscharacters}, and that the upper bound of \eqref{eq: conj Dousse Feray} holds for $k=O(n^{1/2})$ in \cite[Theorem 4]{DousseFeray2019skewdiagramscharacters}. In order to prove the upper bound of \eqref{eq: conj Dousse Feray} via the strategy of Dousse and Féray, one would need to prove the conjecture of Moore and Russell.

\medskip

The probabilistic point of view on the quantity $A_{\lambda,\mu}$ leads to the following elementary improvement of \cite[Theorem 5]{DousseFeray2019skewdiagramscharacters}: the assumption that $\lambda$ is balanced can be removed, and the exponent is improved by a logarithmic factor.

\begin{proposition}
    Let $C>0$. As $n\to \infty$, uniformly over all $\lambda\vdash n$ and all $C$-balanced diagrams $\mu \subset \lambda$, we have 
    \begin{equation}
        A_{\lambda, \mu} \leq (C^2e)^k.
    \end{equation}
\end{proposition}
\begin{proof}
    Let $n$ be a large integer, let $\lambda \vdash n$, and let $\mu\subset \lambda$ be a $C$-balanced Young diagram. By the hook length formula we have 
    \begin{equation}
        d_\mu = \frac{k!}{H(\mu, \mu)} \geq \frac{k!}{(C\sqrt{k})^k}.
    \end{equation}
Moreover $k!\geq (k/e)^k$ and $P_\lambda(\mu) \leq 1$ since $P_\lambda(\cdot)$ is a probability. We conclude that
\begin{equation}
  A_{\lambda, \mu} = \frac{P_\lambda(\mu)}{\Pl_k(\mu)} \leq  \frac{1}{\Pl_k(\mu)} = \frac{d_\mu^2}{k!} \leq \frac{\pg (C\sqrt{k})^k\pd^2}{k!} \leq (C^2e)^k. \qedhere
\end{equation}
\end{proof}

\section{Bounds on Kronecker coefficients via characters}\label{s: Kronecker}

Let $n\geq 2$. Let $X$ be a uniform permutation on $\kS_n$.
The Kronecker coefficient for the representations $\lambda, \mu, \nu \vdash n$ is given by
\begin{equation}
    g(\lambda, \mu, \nu) = \frac{1}{n!}\sum_{\sigma \in \kS_n} \ch^\lambda(\sigma) \ch^\mu(\sigma) \ch^\nu(\sigma) = \bbE\cg \ch^\lambda(X) \ch^\mu(X) \ch^\nu(X)\cd.
\end{equation}
We refer to \cite{PakPanova2023DurfeeKronecker} for more details and history on Kronecker coefficients.

\medskip

Using Theorem \ref{thm: character bound useful for fixed point free intro} (proved in Section \ref{s: Character bound in function of the diagonal length and the number of cycles}), we obtain the following simple bounds on Kronecker coefficients. Recall that $\delta(\lambda)$ and $\ell(\lambda)$ respectively denote the diagonal length and the number of rows of a diagram $\lambda$.

\begin{proposition}\label{prop: bounds on Kronecker coefficients general}
    Let $n\geq 2$. Let $\lambda, \mu, \nu \vdash n$. We have the following.
\begin{enumerate}
    \item $g(\lambda, \mu, \nu) \leq \binom{n+K-1}{n}$, where $K = \ell(\lambda)\ell(\mu) \ell(\nu)$.
    \item  $g(\lambda, \mu, \nu) \leq \binom{n+K-1}{n}$, where $K = 8\delta(\lambda)\delta(\mu) \delta(\nu)$.
\end{enumerate}
\end{proposition}
\begin{proof}
By the triangle inequality and Theorem \ref{thm: character bound useful for fixed point free intro}, for both (a) and (b) we have
\begin{equation}
    g(\lambda, \mu, \nu) \leq \bbE\cg | \ch^\lambda(X)||\ch^\mu(X)| |\ch^\nu(X)|\cd\leq \bbE[K^{\cyc(X)}],
\end{equation}
where $X$ is a uniform permutation of $\kS_n$.
Since the generating function of $\cyc(X)$ is given, for $z\in \bbC$, by $\bbE[z^{{\cyc(X)}}] = \frac{z(z+1)\cdots(z+n-1)}{n!}$, which is equal to $\binom{n+z-1}{n}$ if $z$ is a positive integer, we conclude that
\begin{equation}
   g(\lambda, \mu, \nu) \leq \bbE[K^{\cyc(X)}] = \binom{n+K-1}{n}. \qedhere
\end{equation}
\end{proof}
We note that Proposition \ref{prop: bounds on Kronecker coefficients general} gives essentially the same bound as \cite[Theorem 1.1]{PakPanova2020KroneckerContingency}, and also implies that if $K =\ell(\lambda) \ell(\mu)\ell(\nu) \leq n$ we have $g(\lambda, \mu, \nu) \leq 4^n$, since $\binom{n+K-1}{n}\leq \binom{2n-1}{n}\leq 4^n$.

\medskip

The following corollary of Proposition \ref{prop: bounds on Kronecker coefficients general} recovers that Kronecker coefficients for diagrams with a bounded diagonal length are polynomial. Our constant 8 in the exponent (which we did not try to optimize) is weaker than the constant 4 that of \cite[Theorem 1.4]{PakPanova2023DurfeeKronecker}, but our proof is significantly simpler.

\begin{corollary}
    Let $n\geq 4$. Let $\lambda, \mu, \nu\vdash n$ and set $k = \max(\delta(\lambda), \delta(\mu),\delta(\nu))$. Then
\begin{equation}
    g(\lambda, \mu, \nu) \leq n^{8k^3}.
\end{equation}
\end{corollary}
\begin{proof}
    Set $K = 8k^3$. By Proposition \ref{prop: bounds on Kronecker coefficients general}, we have
\begin{equation}
    g(\lambda, \mu, \nu) \leq \binom{n+K-1}{n} \leq \binom{n+K}{n}.
\end{equation}
Now if $K\geq n$ we have $\binom{n+K}{n}\leq \binom{2K}{n}\leq 2^{2K}= 4^K \leq n^K$, and if $K\leq n$, using that $K\geq 4$ so that $K!\geq 2^K$, we have
\begin{equation}
    \binom{n+K}{n} = \binom{n+K}{K}\leq \binom{2n}{K}\leq \frac{(2n)^K}{K!} \leq n^K. \qedhere
\end{equation}
\end{proof}

\section*{Acknowledgements}

We thank Valentin Féray and Greta Panova for stimulating conversations, and the referees for their comments. The author was supported by the Pacific Institute for the Mathematical Sciences, the Centre national de la recherche scientifique, and the Simons foundation, via a PIMS-CNRS-Simons postdoctoral fellowship.

{\footnotesize
\bibliographystyle{alpha}
\bibliography{bibliographieLucas}
}
\end{document}